\documentclass{scrarticle}

\usepackage[colorlinks=true,linkcolor=teal,citecolor=teal,hypertexnames=false]{hyperref}
\usepackage[margin=2.5cm]{geometry}
\usepackage{amsmath,amssymb,amsthm}
\usepackage{enumitem}
\usepackage{comment}
\usepackage{microtype}
\usepackage[english]{babel}
\usepackage[utf8]{inputenc}
\usepackage{mathrsfs}
\usepackage{bbm}

\usepackage{tikz}
\usetikzlibrary{matrix,arrows}
\usepackage{tikz-cd}
\usetikzlibrary{cd}
\usepackage{comment}

\mathcode`A="7041 \mathcode`B="7042 \mathcode`C="7043 \mathcode`D="7044
\mathcode`E="7045 \mathcode`F="7046 \mathcode`G="7047 \mathcode`H="7048
\mathcode`I="7049 \mathcode`J="704A \mathcode`K="704B \mathcode`L="704C
\mathcode`M="704D \mathcode`N="704E \mathcode`O="704F \mathcode`P="7050
\mathcode`Q="7051 \mathcode`R="7052 \mathcode`S="7053 \mathcode`T="7054
\mathcode`U="7055 \mathcode`V="7056 \mathcode`W="7057 \mathcode`X="7058
\mathcode`Y="7059 \mathcode`Z="705A

\newtheorem{Theorem}{Theorem}[section]
\newtheorem*{Theorem*}{Theorem}
\newtheorem{Lemma}[Theorem]{Lemma}

\newtheorem{Proposition}[Theorem]{Proposition}
\newtheorem*{Proposition*}{Proposition}
\newtheorem{Corollary}[Theorem]{Corollary}
\theoremstyle{definition}
\newtheorem{Definition}[Theorem]{Definition}
\newtheorem{Construction}[Theorem]{Construction}
\newtheorem*{Definition*}{Definition}
\newtheorem{Remark}[Theorem]{Remark}
\newtheorem{Notation}[Theorem]{Notation}

\newtheorem*{rem*}{Remark}

\usepackage{thmtools}
\usepackage{thm-restate}

\newcommand{\AMot}{A\mathrm{Mot}} 

\newcommand{\colim}{\mathrm{colim}}

\newcommand{\Cl}{\mathrm{Cl}}

\newcommand{\reg}{\mathrm{reg}}
\newcommand{\eff}{\mathrm{eff}}

\newcommand{\res}{\operatorname{res}}

\newcommand{\Hom}{\mathrm{Hom}}

\newcommand{\Ext}{\mathrm{Ext}}
\newcommand{\End}{\mathrm{End}}

\newcommand{\id}{\mathrm{id}}
\newcommand{\Snk}{\mathrm{Snk}}

\newcommand{\Gal}{\mathrm{Gal}}

\newcommand{\bF}{\mathbb{F}}

\newcommand{\bQ}{\mathbb{Q}}
\newcommand{\bR}{\mathbb{R}}

\newcommand{\bZ}{\mathbb{Z}}

\usepackage{calrsfs}
\DeclareMathAlphabet{\pazocal}{OMS}{zplm}{m}{n}

\newcommand{\cO}{\pazocal{O}}

\newcommand{\sA}{\mathcal{A}}
\newcommand{\sB}{\mathcal{B}}

\newcommand{\sD}{\mathcal{D}}

\newcommand{\KI}{K_{\infty}}
\newcommand{\CI}{\mathbb{C}_{\infty}}

\newcommand{\Lie}{\mathrm{Lie}}
\newcommand{\FI}{F_{\infty}}

\newcommand{\fj}{\mathfrak{j}}

\newcommand{\fm}{\mathfrak{m}}

\makeatletter
\newsavebox{\@brx}
\newcommand{\llangle}[1][]{\savebox{\@brx}{\(\m@th{#1\langle}\)}%
  \mathopen{\copy\@brx\kern-0.5\wd\@brx\usebox{\@brx}}}
\newcommand{\rrangle}[1][]{\savebox{\@brx}{\(\m@th{#1\rangle}\)}%
  \mathclose{\copy\@brx\kern-0.5\wd\@brx\usebox{\@brx}}}
\makeatother

\numberwithin{equation}{section}

\title{A Barsotti--Weil isomorphism for Taelman units}
\author{Quentin Gazda\footnote{Sorbonne Universit\'e and Universit\'e Paris Cit\'e, CNRS, IMJ-PRG, 75005 Paris, \emph{Email:} \texttt{quentin@gazda.fr}}}
\date{}

\begin{document}

\maketitle

\begin{center}\textsc{Abstract}\end{center}
The positive-characteristic analogue of motivic cohomology with coefficients in an Anderson $A$-motive $\underline{M}$ over a function field of transcendence degree one was introduced in \cite{gazda2}. In this paper, we prove that, when $\underline{M}$ arises as the dual of the $A$-motive of an Anderson module $E$ that is abelian, uniformizable, and has everywhere good reduction, the corresponding $A$-motivic cohomology coincides with the module of Taelman units of $E$, up to a twist by the module of K\"ahler differentials of $A$. To this end, we recast the theory of Anderson generating functions in the spirit of earlier work of Maurischat and the author on special functions.

\tableofcontents

\section{Introduction}
Let $\mathcal{A}$ be an abelian variety over the field of rational numbers $\mathbb{Q}$. The Barsotti--Weil formula implies a canonical isomorphism between the group of $\mathbb{Q}$-valued points of $\mathcal{A}$ and the extension group of commutative fppf sheaves on the big fppf site over $\bQ$ between the dual abelian variety of $\mathcal{A}$ and $\mathbb{G}_m$:
\[
\mathcal{A}(\mathbb{Q}) \cong \Ext^1_{\mathrm{AbGrpSch}/\mathbb{Q}} \bigl(\mathcal{A}^{\vee},\mathbb{G}_m\bigr);
\]
\emph{e.g.} \cite[Lemma~1.13.3]{barbieri-viale-kahn2016}. Its motivic counterpart is the expected canonical isomorphism
\begin{equation}\label{eq:number-field-BW}
\mathcal{A}(\mathbb{Q})\otimes_{\bZ}\bQ \cong \Ext^1_{\mathrm{MxdMot}_{\mathbb{Q}}}(\mathbbm{1},h^1(\mathcal{A})^{\vee})
\end{equation}
where the right-hand side stands for the extension group in a hypothetical $\mathbb{Q}$-linear category of mixed motives. There is, for instance, an unconditional $1$-motive incarnation of the above (see \cite[\S 10]{deligne1974} and \cite[Theorem~2.1.2 and Appendix~C, Proposition~C.8.3(b)]{barbieri-viale-kahn2016}). By the Mordell--Weil theorem, the left-hand side of \eqref{eq:number-field-BW} is a finite-dimensional $\bQ$-vector space, which makes this interpretation particularly useful in arithmetic applications.

The aim of this paper is to prove two Barsotti--Weil-style isomorphisms in function field arithmetic, where abelian varieties are replaced by Anderson modules. Let $\bF$ be a finite field with $q$ elements and let $(C,\cO_C)$ be a geometrically irreducible smooth projective curve over $\bF$ marked with a closed point $\infty$. Consider the following rings of algebraic functions on $C$:
\[
A:=\cO_C(C\setminus\{\infty\}) \quad \text{and} \quad K:=\operatorname{Frac}(A)=\mathbb F(C).
\]
Anderson introduced analogues of abelian varieties with coefficients in $A$, as higher-dimensional generalizations of Drinfeld modules \cite{anderson}. An Anderson module $E$ over an $A$-algebra $R$ is a particular instance of a smooth $A$-module scheme over $R$ (see Subsection~\ref{subsec:recollection} for a complete definition). Anderson developed their theory of motives, called \emph{$A$-motives}, which consist of certain $A\otimes_{\bF} R$-modules $M$ equipped with a semilinear endomorphism of $M$ with respect to the $q$-Frobenius on $R$. If $E$ is an Anderson module over an $A$-algebra $R$, its associated $A\otimes R[\tau]$-module is
\[
M(E):=\Hom_{\bF}(E,\mathbb{G}_a),
\]
where homomorphisms are those of $\bF$-vector space schemes over $R$. The functor $M$ shares many similarities with the functor $h^1$ in the number field setting.

We establish in this paper a first Barsotti--Weil-type isomorphism (see Theorem \ref{thm:barsotti-weil}).
\begin{Theorem*}
Let $E$ be an Anderson module over an $A$-algebra $R$ and let $M(E)$ denote its associated $A\otimes R[\tau]$-module. Then we have
\[
E(R) \cong \Ext^{1}_{A\otimes R[\tau]}(M(E),A\otimes R)\otimes_A \Omega_{A/\bF}^1
\]
where the extension module is taken in the category of $A\otimes_{\bF} R[\tau]$-modules. When $E$ is abelian, this implies
\begin{equation}\label{eq:function-field-BW}
E(R) \cong \Ext^{1,\mathrm{reg}}_{\AMot_R}(\mathbbm{1},\underline{M}(E)^{\vee})\otimes_A \Omega_{A/\bF}^1
\end{equation}
where the extension module is now taken in the category of $A$-motives over $R$; the superscript ``reg'' means \emph{regulated extensions}, a notion introduced in \cite{gazda} to allow the construction of a regulator map.\footnote{\label{fn:regulation}Regulation is originally defined over fields in \cite[Definition~5.1]{gazda}. For a general $A$-algebra $R$, we use the criterion of \cite[Proposition~5.3]{gazda} as a definition. Given an extension $0\to\underline{M}\to\underline{X}\to\underline{N}\to0$ of $A$-motives over $R$, choose an $A\otimes R$-linear splitting $X\simeq M\oplus N$, which exists because $N$ is projective. Let $v:\tau^*N[\mathfrak{j}^{-1}]\to M[\mathfrak{j}^{-1}]$ be the upper diagonal entry of $\tau_X$ in this splitting. The extension is called \emph{regulated} if $v=\tau_M\circ f-g\circ\tau_N$ after inverting $\mathfrak{j}$, for $A\otimes R$-linear maps $f:\tau^*N\to\tau^*M$ and $g:N\to M$. This condition is independent of the splitting: changing it by $h:N\to M$ replaces $v$ by $v+\tau_M\circ\tau^*h-h\circ\tau_N$, which preserves the condition.}
\end{Theorem*}
This theorem was known in several cases. If $E$ is abelian, $A=\bF[t]$ and $R$ is an algebraically closed field, it was already established by Taelman \cite{1-t-motif}. It was also established for Drinfeld modules by Mornev \cite{mornev-shtuka} (in the language of shtukas), and for abelian Anderson modules over perfect $R$ by the author and Maurischat \cite{gazda-maurischat-pairing}.

Let $F$ be a finite extension of $K$ and let $\mathcal{O}_F$ be the integral closure of $A$ in $F$. In contrast with the number-field situation, the finite generation of $E(\mathcal{O}_F)$ fails already in dimension one. The correct replacement is the submodule of $E(\mathcal{O}_F)$ consisting of \emph{Taelman units}, denoted by $U(E/\mathcal{O}_F)$ (see \cite{taelman-dirichlet} for Drinfeld modules). A natural problem is then to identify the extension submodule to which it corresponds under the isomorphism \eqref{eq:function-field-BW}. The aim of this paper is to show that $U(E/\mathcal{O}_F)$ agrees with the submodule of extensions having analytic reduction at $\infty$, introduced in \cite{gazda2}.

In order to formulate our statement precisely, let us introduce additional notation. We write $\KI$ for the completion of $K$ at $\infty$, equipped with $v_{\infty}$ and let $F_{\infty}$ be $F\otimes_K \KI$. After fixing a separable closure $\KI^s$ of $\KI$, we also write $F_{\infty}^s$ for $F\otimes_K \KI^s$. Observe that since $\KI$ is a separable field extension of $K$, both $F_{\infty}$ and $F_{\infty}^s$ are finite products of fields.

We now assume that $E$ is a uniformizable Anderson $A$-module over $\cO_F$; by definition, this means that its exponential map $\exp_E:\Lie_E(\FI^s)\to E(\FI^s)$ is surjective. This map is equivariant with respect to the profinite group $G_{\infty}:=\Gal(K_{\infty}^s|K_{\infty})$. Taking invariants under its action and using additive Hilbert 90 yields a surjective connecting homomorphism to continuous Galois cohomology
\begin{equation}\label{eq:connecting}
E(F_{\infty}) \longrightarrow H^1(G_{\infty},\Lambda_E)
\end{equation}
where $\Lambda_E:=\ker(\exp_E)$ is the period lattice of $E$, see \cite[\S 2]{taelman-dirichlet}. We define Taelman's map $t_{\exp}$ as the restriction of \eqref{eq:connecting} to $E(\cO_F)$, and set
\[
U(E/\cO_F):=\ker(t_{\exp}) \quad \text{and} \quad H(E/\cO_F):=\mathrm{coker}(t_{\exp}),
\]
respectively for the \emph{Taelman unit module} and the \emph{Taelman class module} of $E$. Taelman's unit theorem \cite[Theorem~1]{taelman-dirichlet} can be stated as follows.
\begin{Theorem}[Taelman]
As $A$-modules, $U(E/\cO_F)$ is finitely generated and $H(E/\cO_F)$ is finite.
\end{Theorem}

The second Barsotti--Weil-type isomorphism connects Taelman's modules to the motivic invariants introduced in \cite{gazda2}. For this, let us further assume that $E$ is abelian over $\mathcal{O}_F$; in particular, $E_F$ has everywhere good reduction. By definition, this means that $\underline{M}(E)$ is a finite projective $A\otimes \cO_F$-module. Because $E$ is also assumed to be uniformizable, $\underline{M}(E)_F$ is rigid analytically trivial (\emph{cf}. \cite{anderson,bockle-hartl}). To lighten notation, we set $\underline{M}:=\underline{M}(E)_F^{\vee}$ which is also rigid analytically trivial with good reduction. The ``Betti realization'' functor from the category of uniformizable $A$-motives to that of continuous $G_{\infty}$\nobreakdash-representations is exact, and therefore induces a map
\[
r_B:\Ext^{1,\reg}_{\cO_F}(\mathbbm{1},\underline{M}) \longrightarrow H^1(G_{\infty},\underline{M}_B)
\]
when restricted to \emph{regulated and integral extensions} (see \cite{gazda} for definitions and \cite[(3.6)]{gazda2} for details on the map $r_B$). Our main theorem is the following.
\begin{Theorem}\label{thm:comparison-taelman}
There is a natural commutative square of $A$-modules
\[
\begin{tikzcd}[column sep=3em]
\Ext^{1,\reg}_{\cO_F}(\mathbbm{1},\underline{M})\otimes_A \Omega_{A/\bF}^1 \arrow[r,"r_B\otimes\mathrm{id}_{\Omega}"]\arrow[d,"\wr"] &  H^1(G_{\infty},\underline{M}_B)\otimes_A \Omega_{A/\bF}^1 \arrow[d,"\wr"] \\
E(\cO_F) \arrow[r,"t_{\exp}"] & H^1(G_{\infty},\Lambda_E) 
\end{tikzcd}
\]
whose vertical arrows are isomorphisms, as indicated. 
\end{Theorem}

Writing $\Ext^{1,\reg}_{\cO_F,\infty}(\mathbbm{1},\underline{M}):=\ker(r_B)$ and $\Cl(\underline{M}):=\mathrm{coker}(r_B)$, the theorem implies that $\Ext^{1,\reg}_{\cO_F,\infty}(\mathbbm{1},\underline{M})\otimes_A\Omega_{A/\bF}^1$ is naturally isomorphic to Taelman's unit module $U(E/\cO_F)$ and that $\Cl(\underline{M})\otimes_A\Omega_{A/\bF}^1$ is naturally isomorphic to $H(E/\cO_F)$. Indeed, tensoring with the invertible $A$-module $\Omega_{A/\bF}^1$ preserves kernels and cokernels. In particular, \cite[Theorem~6.1]{gazda2} recovers Taelman's theorem.

This shows that the formalism of $A$-motivic cohomology provides a reasonable definition of Taelman units/class modules for $A$-motives that do not arise from Anderson $A$-modules, or for abelian $A$-modules not having everywhere good reduction. An important benefit is that $A$-motivic cohomology might be easier to compute than Taelman units: in the case of Carlitz tensor powers where $A=\bF[t]$, the former has been partly computed in \cite{gazda-maurischat-ext}, and the following corollary can be deduced.

\begin{Corollary}\label{cor:carlitz-symbols}
Assume that $A=\bF[t]$. Let $k_F$ be the constant field of $F$, and let $g$ be its genus over $k_F$. For each place $v$ of $F$ above $\infty$, let $e_v$ be its ramification index and let $\deg_{k_F}(v)$ be its residue degree over $k_F$. Let $n\geq 1$ satisfy
\[
\sum_{v\mid\infty}\left(\frac{qn}{q-1}e_v-1\right)\deg_{k_F}(v)>2g-2,
\]
and let $C^{\otimes n}$ denote the $n$-th tensor power of the Carlitz module over $\cO_F$. For $\alpha\in \cO_F$, denote by $[\alpha]_n:=(0,\ldots,0,-\alpha)$ the \emph{symbol of $\alpha$} in $C^{\otimes n}(\cO_F)$. If $v(\alpha)>-\frac{nq}{q-1}$ for all valuations $v$ of $F$ above $\infty$, normalized by $v|_K=v_\infty$, then $[\alpha]_n$ is a Taelman unit and $U(C^{\otimes n}/\cO_F)$ is generated by such symbols over $\bF[t]$. 
\end{Corollary}

The strategy of the proof of Theorem \ref{thm:comparison-taelman} is to write down the exponential map of $E$ in terms of its motive $\underline{M}(E)$. This idea goes back to Sinha's appendix \cite{sinha} and was later developed by Taelman \cite{1-t-motif} using an alternative approach. Our construction is different and is rooted in the abstract theory of Anderson generating functions developed here. As another distinction, we work with the general coefficient ring $A$, where the module of K\"ahler differentials plays a decisive role. 

After a brief review of Anderson modules and homological preliminaries, we give the proof of the first Barsotti--Weil-type isomorphism in Section \ref{sec:generalities}. In Section \ref{sec:anderson-generating-functions}, we develop the theory of Anderson generating functions, extending previous work of the author and Maurischat \cite{gazda-maurischat} on special functions. Given an Anderson module $E$ over $\FI$ and a Banach $\FI$-algebra $B$, a generating function for its base change to $B$ is an element $\psi\in E(B\langle A\rangle)$ that satisfies 
\[
\text{for~all~}a\in A:\quad (1\otimes \phi_a)\psi-(a\otimes 1)\psi \in A\otimes E(B).
\]
We show that, up to a differential form, each generating function can be described as the generating series of a division tower, as considered by Anderson. The exponential map is recovered by comparing the residues of generating functions at the characteristic graph and at $\infty$. Section \ref{sec:proof} is devoted to the proof of Theorem~ \ref{thm:comparison-taelman} and Corollary~\ref{cor:carlitz-symbols}. We conclude in Section \ref{sec:regulators} by proving that the Hodge-Pink regulator introduced in \cite{gazda2} generalizes Taelman's regulator. 

\paragraph{Acknowledgments.}
The idea of seeking a motivic incarnation of the exponential map was suggested to us by Tuan Ngo Dac. The author is grateful for this suggestion and hopes that it will find further applications in the future. Parts of this work were presented by the author during the thematic programme \emph{Arithmetic of Global Function Fields}, held at ICMAT in June 2026. The author thanks the organizers and participants for their comments during the lectures. This work was funded by the project ANR-25-CE40-3307-01, \emph{Al$K$traZ: Algebraic K-Theory, Trace Maps and the Zagier Conjecture}.

\section{Generalities on Anderson modules}\label{sec:generalities}
\subsection{Recollection}\label{subsec:recollection}
Let $R$ be an $A$-algebra via a structure morphism $\ell:A\to R$. For our purposes, \emph{an Anderson $A$-module $E$ over $R$} is a smooth and affine $A$-module scheme over $R$ satisfying the following two properties:
\begin{enumerate}[label=$(\mathrm{And})_{\arabic*}$]
\item\label{item:zariski-form} As an $\bF$-module scheme over $R$, $E$ is Zariski-locally isomorphic to $\mathbb{G}_a^d$ for some $d>0$;
\item\label{item:tangentially-compatible} For all $a\in A$, $\partial \phi_a-\ell(a)$ acts nilpotently on $\Lie_E(R)$, where $\phi_a$ denotes the action of $a$ on $E$ and $\partial \phi_a$ its induced action on $\Lie_E(R)$. 
\end{enumerate}

To $E$, we associate its \emph{motive} $M(E):=\Hom_{\bF}(E,\mathbb{G}_a)$; \emph{i.e.} the group of homomorphisms of $\bF$-module schemes over $R$ from $E$ to $\mathbb{G}_a$. It acquires an action of $A$ by acting on $E$ and an action of $R$ by acting on $\End(\mathbb{G}_a)$; both actions coincide on $\bF$, making $M(E)$ an $A\otimes R$-module\footnote{Throughout, unlabelled tensor products are over $\mathbb{F}$.}. 

Because $\End_{\bF}(\mathbb{G}_a)=R[\tau]$, where $\tau$ corresponds to the $q$-power Frobenius endomorphism of $\mathbb{G}_a$, the hypothesis \ref{item:zariski-form} makes $M(E)$ a finitely generated projective\footnote{Note that this is false if one merely assumes instead of \ref{item:zariski-form} that $E$ is a flat form of $\mathbb{G}_a^d$, as in \cite{hartl-isogeny}.} left $R[\tau]$-module of constant rank $d$. The action of $\tau$ can also be made linear, giving rise to an $A\otimes R$-linear map
\[
\tau_E:\tau^* M(E)\longrightarrow M(E)
\]
whose cokernel, thanks to \ref{item:tangentially-compatible}, is killed by a power of the ideal $\mathfrak{j}:=\ker(A\otimes R\to R)$. 
\begin{Definition}[Abelian $A$-modules]
We call $E$ \emph{abelian} if $M(E)$ is finite projective over $A\otimes R$ of constant rank; equivalently, if the pair $\underline{M}(E):=(M(E),\tau_E)$ defines an $A$-motive over $R$ in the sense of \cite[Definition~3.2]{gazda2}. 
\end{Definition}

The data of $\underline{M}(E)$ is sufficient to recover $E$. If $S$ is an $R$-algebra, we have a canonical $A$-linear isomorphism \cite[Proposition~3.2(1)]{gazda-maurischat-pairing}
\begin{equation}\label{eq:fromM(E)toE}
E(S) \cong \Hom_{R[\tau]}(M(E),S)
\end{equation}
where $S$ is treated as a left $R[\tau]$-module via $\tau\cdot s:=s^q$. The isomorphism \eqref{eq:fromM(E)toE} maps $e_S\in E(S)$ to $(m\in M(E)\mapsto m(e_S)\in S)$.

\begin{Notation}[Extended functor of points]\label{not:extended-functor-of-points}
Formula \eqref{eq:fromM(E)toE} allows us to consider points of $E$ not only over $R$-algebras, but also over left $R[\tau]$-modules $T$ by setting
\begin{equation}\label{eq:functor-of-points}
E(T):=\Hom_{R[\tau]}(M(E),T).
\end{equation}
This perspective was first used in \cite{caruso-gazda-lucas} and we will use it constantly throughout. If $T$ is a left $R[\tau]$-module on which $\tau$ acts by zero, we shall write
\[
\Lie_E(T):=E(T).
\]
This is consistent with the usual functor $\Lie_E$. When we write $\Lie_E(T)$ for an $R$-module $T$, it is tacitly assumed that we treat $T$ as a left $R[\tau]$-module where $\tau$ acts by zero.
\end{Notation}

We end this subsection with the following lemma, which is well-known in the abelian case, \emph{e.g.} \cite[Corollary~6.6]{gazda-maurischat-pairing}.
\begin{Lemma}\label{lem:no-non-zero-morphism}
Let $M$ be an $A\otimes R[\tau]$-module which is finitely generated over $R[\tau]$. Let $u\in A$ be nonconstant. Then $\Hom_{R[u,\tau]}(M,A\otimes R)=(0)$. In particular, $\Hom_{A\otimes R[\tau]}(M,A\otimes R)=(0)$.
\end{Lemma}
\begin{proof}
First, we equip $A$ with the degree $\deg=-v_{\infty}$. Given a finite-dimensional $\bF$-subspace $V$ of $R$, we extend $\deg$ to $A\otimes V$ in the classical way; \emph{i.e.} by taking the coordinate-wise maximum of the degree with respect to a chosen basis of $V$ over $\bF$ (although it does not depend on such a choice). Then, we extend $\deg$ to $A\otimes R$ by writing $R$ as $\colim_{V\subset R} V$, indexed over finite-dimensional $\bF$-vector spaces inside $R$. Clearly, we have
\begin{align*}
&\text{for }a\in A,~h\in A\otimes R, &&\deg(a\cdot h)=\deg(a)+\deg(h), \\
&\text{for }h_1,h_2\in A\otimes R, &&\deg(h_1+h_2)\leq \max\{\deg(h_1),\deg(h_2)\}.
\end{align*}
Let $f\in \Hom_{R[u,\tau]}(M,A\otimes R)$. We fix a finite family $(m_1,\ldots,m_g)$ which generates $M$ over $R[\tau]$ and set $h_i:=f(m_i)\in A\otimes R$, for $i\in\{1,\ldots,g\}$. If $f\neq 0$, let $j\in \{1,\ldots,g\}$ be such that $\deg(h_j)$ is maximal. Because $um_j\in M$, there exist $p_1(\tau),\ldots,p_g(\tau)\in R[\tau]$ such that
\[
um_j=\sum_{i=1}^g{p_i(\tau)m_i} \quad \text{and thus} \quad uh_j=\sum_{i=1}^g{p_i(\tau)h_i}
\]
But the left-hand side of the latter equality has degree $>\deg(h_j)$ while the right-hand side has degree $\leq \deg(h_j)$, a contradiction. Therefore $f=0$.
\end{proof}

\subsection{Functor of $A$-invariants}
We call $u\in A$ \emph{separating} if $K/\mathbb F(u)$ is finite separable; equivalently, $\mathbb F[u]\hookrightarrow A$ is finite flat and generically étale. Given a separating element $u\in A$, we define \emph{the dualizing module of $A/\bF[u]$:}
\[
\sD_{A/\bF[u]}:=\{x\in A\otimes_{\bF[u]} A~|~\text{for~all~}a\in A:~(a\otimes_{\bF[u]} 1-1\otimes_{\bF[u]} a)x=0\}.
\]
It is an $A$-module (with respect to the left or right multiplication by $A$, which coincide by construction). Under the multiplication map $A\otimes_{\bF[u]}A \to A$ the module $\sD_{A/\bF[u]}$ corresponds to the different ideal $\mathfrak{d}_{A/\bF[u]}$. From generalities on K\"ahler differentials (\emph{e.g.} \cite[\S 3.2]{gazda-maurischat}), we have a canonical isomorphism
\begin{equation}\label{eq:multiplication-different}
\sD_{A/\bF[u]}\otimes_A \Omega_{A/\bF}^1\xrightarrow{\sim} Adu, \quad x\otimes_A \omega\longmapsto x\omega.
\end{equation}
For later use, let 
\begin{equation}\label{eq:dualizing-element}
Du:=\sum_{i}{d_i \otimes_A \omega_i}\in \sD_{A/\bF[u]}\otimes_A \Omega_{A/\bF}^1, \quad d_i\in \sD_{A/\bF[u]}, \quad \omega_i\in \Omega_{A/\bF}^1
\end{equation}
be the inverse image of $du$ through \eqref{eq:multiplication-different} (\emph{the dualizing element}).

More generally, we will require the \emph{functor of $A$-invariants}.
\begin{Definition}[Functor of $A$-invariants]\label{def:functor-of-A-invariants}
We call the functor
\begin{align*}
(-)^{A}~: \{(A\otimes_{\bF[u]}A)-\text{modules}\} &\longrightarrow \{A-\text{modules}\},\\ 
 V &\longmapsto V^A:=\{v\in V~|~(a\otimes 1)v=(1\otimes a)v\text{ for all }a\in A\}
\end{align*}
which sends an $A\otimes_{\bF[u]}A$-module to its submodule on which both $A$-actions coincide, the \emph{functor of $A$-invariants}.
\end{Definition}

This functor maps the unit module $A\otimes_{\mathbb{F}[u]}A$ to $\sD_{A/\bF[u]}$ which, as mentioned above, is isomorphic to the dual of $\Omega^1_{A/\mathbb{F}}$. More generally, we have:

\begin{Lemma}\label{lem:A-invariants-D}
Let $T$ be an $A$-module. Then $(A\otimes_{\bF[u]}T)^A =\sD_{A/\bF[u]}\otimes_A T$.
\end{Lemma}

\begin{proof}
Consider the $A$-module
\[
\omega_{A/\bF[u]}:=\Hom_{\bF[u]}(A,\bF[u]).
\]
Since $A$ is finite projective over $\bF[u]$, evaluation induces, for every $A$-module $T$, a canonical isomorphism
\[
A\otimes_{\bF[u]} T \xrightarrow{\sim} \Hom_{\bF[u]}(\omega_{A/\bF[u]},T), \qquad a\otimes t\longmapsto \bigl(\lambda\mapsto\lambda(a)t\bigr).
\]
Under this isomorphism, multiplication by $a\in A$ on the first factor of $A\otimes_{\bF[u]} T$ corresponds to $f(\lambda)\mapsto f(a\lambda)$, whereas multiplication by $a$ on $T$ corresponds to $f(\lambda)\mapsto a f(\lambda)$. Consequently, the relation
\[
\text{for~}v\in A\otimes_{\bF[u]} T: \quad  (a\otimes1)v=(1\otimes a)v \qquad\text{for every }a\in A
\]
is equivalent to $A$-linearity of the corresponding map $\omega_{A/\bF[u]}\to T$. Therefore
\[
(A\otimes_{\bF[u]} T)^A \cong \Hom_A(\omega_{A/\bF[u]},T).
\]
Taking $T=A$ gives
\[
\sD_{A/\bF[u]} =(A\otimes_{\bF[u]}A)^A \cong \Hom_A(\omega_{A/\bF[u]},A).
\]
Moreover, $\omega_{A/\bF[u]}$ is a finite torsion-free $A$-module of rank one. Since $A$ is a Dedekind domain, $\omega_{A/\bF[u]}$ is invertible. Hence
\[
\Hom_A(\omega_{A/\bF[u]},T) \cong \Hom_A(\omega_{A/\bF[u]},A)\otimes_A T \cong \sD_{A/\bF[u]}\otimes_A T.
\]
Combining the preceding isomorphisms yields
\[
(A\otimes_{\bF[u]} T)^A \cong \sD_{A/\bF[u]}\otimes_A T,
\]
as claimed.
\end{proof}

\begin{Proposition}\label{prop:extension-as-A-invariants}
Let $M$ be an $A\otimes R[\tau]$-module, finitely generated over $R[\tau]$. Then 
\[
\Ext^1_{A\otimes R[\tau]}(M,A\otimes R) = \Ext^1_{R[u,\tau]}(M,A\otimes R)^A.
\]
\end{Proposition}
\begin{proof}
We begin with the formula $\Hom_{R[u,\tau]}(M,-)^A=\Hom_{A\otimes R[\tau]}(M,-)$, where both $(-)^A$ and $\Hom_{R[u,\tau]}(M,-)$ are additive, left-exact functors. We write $(-)^{R^i A}$ for the $i$th right derived functor of $(-)^A$. The standard hypothesis are satisfied to apply the Grothendieck spectral sequence theorem, and from the first-quadrant five terms exact sequence we get:
\[
\Hom_{R[u,\tau]}(M,-)^{R^1 A} \to \Ext^1_{A\otimes R[\tau]}(M,-) \xrightarrow{d} \Ext^1_{R[u,\tau]}(M,-)^A \to \Hom_{R[u,\tau]}(M,-)^{R^2 A} \to \cdots
\]
where $d$ is the canonical map. Applied to the left $A\otimes R[\tau]$-module $A\otimes R$ and using Lemma \ref{lem:no-non-zero-morphism}, the first and fourth terms vanish and we obtain that $d$ is an isomorphism. 
\end{proof}

\subsection{Barsotti--Weil isomorphism}
The Barsotti--Weil isomorphism relates the module of points of $E$ to extension modules in categories of $A$-motives. To our knowledge, the strongest form of the statement is originally due to Taelman \cite[Theorem 8.1.1]{1-t-motif}, generalizing prior work of Papanikolas and Ramachandran \cite{PR}, and Woo \cite{woo}. We also refer the reader to the recent work of Głoch--K\c{e}dzierski--Krasoń \cite{gkk,kk,kp} for further results in this direction. The formula was subsequently generalized for general coefficient rings $A$ by Mornev \cite[Theorem 9.4]{mornev-shtuka}. For perfect $A$-algebras, it was obtained as a consequence of Hartl--Juschka's pairing in \cite{gazda-maurischat-pairing}. We prove it in full generality here, without abelianness assumptions.

\begin{Theorem}[Barsotti--Weil isomorphism]\label{thm:barsotti-weil}
Let $E$ be a (not necessarily abelian) Anderson $A$-module over an $A$-algebra $R$. There is an isomorphism of $A$-modules
\[
E(R) \cong \Ext^{1}_{A\otimes R[\tau]}(M(E),A\otimes R)\otimes_A \Omega_{A/\bF}^1.
\]
\end{Theorem}
\begin{proof}
We start with the identification \eqref{eq:fromM(E)toE}. Fix a separating element $u\in A$. Writing $M=M(E)$, we have 
\[
A\otimes E(R) =\Hom_{R[\tau]}(M,A\otimes R),
\]
where $A\otimes R$ is treated as a left $R[\tau]$-module with the $A$-linear action of $\tau$. Hence
\begin{align*}
A\otimes_{\mathbb{F}[u]}E(R) &= \frac{A\otimes E(R)}{(u\otimes 1-1\otimes \phi_u)(A\otimes E(R))} \\
&= \frac{\operatorname{Hom}_{R[\tau]}(M,A\otimes R)}{\left\{(u\otimes 1)g(\ast)-g\bigl((1\otimes u)\ast\bigr)~\middle|~g\in \operatorname{Hom}_{R[\tau]}(M,A\otimes R)\right\}}
\end{align*}
where the second isomorphism follows from the $A$-linearity of the identification \eqref{eq:fromM(E)toE}. Using that $M$ is projective over $R[\tau]$, we see that the right-hand side computes extensions of $M$ by $A\otimes R$ in the category of $R[u,\tau]$-modules (the underlying sequence of $R[\tau]$-modules splits, so the action of $u$ is parametrized by an $R[\tau]$-linear map from $M$ to $A\otimes R$). We obtain an isomorphism
\begin{equation}\label{eq:BW-for-u}
A\otimes_{\mathbb{F}[u]}E(R)
\simeq
\operatorname{Ext}^1_{R[u,\tau]}(M,A\otimes R).
\end{equation}
Applying $(-)^{A}$ to \eqref{eq:BW-for-u} and appealing to Lemma \ref{lem:A-invariants-D} and Proposition \ref{prop:extension-as-A-invariants}, we get
\[
\sD_{A/\bF[u]}\otimes_A E(R)
\cong 
\operatorname{Ext}^1_{A\otimes R[\tau]}(M,A\otimes R),
\]
which is the desired Barsotti--Weil formula after twisting by $\Omega^1_{A/\mathbb{F}}$.
\end{proof}

\begin{Corollary}
Let $E$ be an abelian Anderson $A$-module over $R$ and let $\underline{M}(E)^{\vee}$ be the dual of its $A$-motive. Then, there is an isomorphism
\[
E(R)\cong \Ext^{1,\reg}_{\AMot_R}(\mathbbm{1},\underline{M}(E)^{\vee})\otimes_A \Omega_{A/\bF}^1.
\]
\end{Corollary}
\begin{proof}
The category of effective $A$-motives embeds fully and faithfully into the category of $A\otimes R[\tau]$-modules and its image is closed under extensions. In particular,
\[
\Ext^1_{A\otimes R[\tau]}(M(E),A\otimes R)\cong \Ext^{1,\eff}_{\AMot_R}(\underline{M}(E),\mathbbm{1}),
\]
where the superscript ``eff'' means extensions of effective $A$-motives. It follows directly from the definition in footnote~\ref{fn:regulation} that an extension of $\underline{M}(E)$ by $\mathbbm{1}$ in $\AMot_R$ is effective if and only if it is regulated. Consequently, using \cite[Corollary~5.4]{gazda}, we get
\[
\Ext^{1,\eff}_{\AMot_R}(\underline{M}(E),\mathbbm{1})\cong \Ext^{1,\reg}_{\AMot_R}(\underline{M}(E),\mathbbm{1})\cong \Ext^{1,\reg}_{\AMot_R}(\mathbbm{1},\underline{M}(E)^{\vee}). 
\]
The desired statement therefore follows from Theorem \ref{thm:barsotti-weil}.
\end{proof}

\begin{Remark}[Computing the Barsotti--Weil isomorphism]\label{rem:compute-barsotti-weil}
Assume that $E$ is abelian. In \cite[Proposition~5.3]{gazda}, the module of effective/regulated extensions was rather presented as a quotient of $\Hom_{A\otimes R}(\tau^* M(E),A\otimes R)$ via a surjective map
\[
\iota: \Hom_{A\otimes R}(\tau^* M(E),A\otimes R)\twoheadrightarrow \Ext^1_{A\otimes R[\tau]}(M(E),A\otimes R)
\]
assigning to $f$ the class of the extension $(A\otimes R)\oplus M(E)$ (which splits at the level of $A\otimes R$-modules) where the action of $\tau$ is represented by a $2$-by-$2$ upper triangular matrix 
\[
\begin{pmatrix} \tau_{\mathbbm{1}} & f \\ 0 & \tau_M \end{pmatrix}.
\]
Here, we describe how to compute the Barsotti--Weil map under this presentation, in the case where $A=\bF[t]$ (and $u$ is chosen to be $t$).

Let $e$ be an element in $E(R)$, seen as a morphism in $\Hom_{R[\tau]}(M(E),R)$. We deduce from the proof of Theorem \ref{thm:barsotti-weil} that the Barsotti--Weil isomorphism maps $e$ to the extension $H=H_e$ of $M(E)$ by $R[t]$ given as follows: $H=R[t]\oplus M(E)$, where the action of $\tau$, denoted by $\tau_H$, is diagonal and where the action of $t$, denoted by $t_H$, is given by a $2$-by-$2$ upper triangular matrix
\[
t_H=\begin{pmatrix} t & e \\ 0 & t_{M(E)} \end{pmatrix}.
\]
To describe $H_e$ as arising from the image of $\iota$, select an $R$-linear map $h:M(E)\to R[t]$ such that
\[
e=t\circ h-h\circ t_M.
\]
Consider the $R$-module $H_e':=R[t]\oplus M(E)$ equipped with the following $R[t,\tau]$-module structure: $t_{H'}$ acts diagonally but where $\tau_{H'}$ is given by the $2$-by-$2$ upper triangular matrix 
\[
\tau_{H'}=\begin{pmatrix} \tau_{\mathbbm{1}} & h\circ \tau_{M(E)}-\tau_{\mathbbm{1}}\circ \tau^*h \\ 0 & \tau_M \end{pmatrix}.
\]
By design, $H_{e}'$ is an extension of $M(E)$ by $R[t]$ as $R[t,\tau]$-modules, and $\left(\begin{smallmatrix} 1 & h \\ 0 & 1 \end{smallmatrix}\right):H_e\to H_e'$ defines an $R[t,\tau]$-linear equivalence of extensions. The equivalence class of $H_{e}'$ is that of $\iota(h\circ \tau_{M(E)}-\tau_{\mathbbm{1}}\circ \tau^*h)$.

In conclusion, the Barsotti--Weil isomorphism maps $e\in \Hom_{R[\tau]}(M(E),R)$ to the image through $\iota$ of 
\[
h\circ \tau_{M(E)}-\tau_{\mathbbm{1}}\circ \tau^*h \in \Hom_{R[t]}(\tau^*M(E),R[t])
\]
where $h:M(E)\to R[t]$ is any $R$-linear map such that $e=t\circ h-h\circ t_{M(E)}$.
\end{Remark}

\begin{Remark}[The Barsotti--Weil isomorphism for Carlitz twists]\label{rem:compute-barsotti-weil-carlitz}
We apply the recipe of Remark \ref{rem:compute-barsotti-weil} to compute the (inverse of the) Barsotti--Weil isomorphism in the case of the Carlitz $n$th twist, where $n\geq 1$. Recall that the latter is the Anderson $t$-module $C^{\otimes n}$ whose $S$-points, for an $R$-algebra $S$, are given by the $\bF$-vector space $S^n$ together with the action of $t$ given by
\[
\text{for~}(x_1,\ldots,x_n)\in S^n:\quad t\cdot (x_1,\ldots,x_n):=(\theta x_1+x_2,\ldots,\theta x_{n-1}+x_n, \theta x_n + x_1^q),
\] 
where $\theta:=\ell(t)$. Given $\alpha\in S$, we set $[\alpha]_n:=(0,\ldots,0,-\alpha)\in C^{\otimes n}(S)$. We call $[\alpha]_n$ the \emph{symbol of $\alpha$}. It can be easily verified directly from the formula that $C^{\otimes n}(S)$ is generated by symbols as an $\bF[t]$-module, although the content of the present remark reproves this.

We denote by $\kappa_i:C^{\otimes n}\to \mathbb{G}_a$ the projection onto the $i$th coordinate. The module $M(C^{\otimes n})$ is free of rank one over $R[t]$, generated by the first projection $\kappa_1$, and $\kappa_i=(t-\theta)^{i-1} \kappa_1$ for $i\in \{1,\ldots,n\}$. In particular, an element $\alpha\in R$ defines the $R[t]$-linear map $f_\alpha\in\Hom_{R[t]}(\tau^*M(C^{\otimes n}),R[t])$ characterized by $f_\alpha(\tau^*\kappa_1)=\alpha$, and therefore an effective/regulated extension $\iota(\alpha)$ of $\underline{M}(C^{\otimes n})$ by $\mathbbm{1}$. The extension $\iota(\alpha)$ was named \emph{the polylog class of $\alpha$} in \cite{gazda-maurischat-ext}. We claim that the Barsotti--Weil isomorphism makes $[\alpha]_n$ correspond to $\iota(\alpha)$.

To prove this claim, consider the unique $R$-linear map $h:M(C^{\otimes n})\to R[t]$ satisfying
\begin{enumerate}
\item for $j\in \{0,\ldots,n-1\}$, $h((t-\theta)^j\kappa_1)=0$,
\item for all $m\in M(C^{\otimes n})$, $h(\tau m)-\tau(h(m))=f_\alpha(\tau^*m)$. 
\end{enumerate}
Since $\kappa_1,\ldots,\kappa_n$ form an $R[\tau]$-basis, these conditions determine $h$ recursively. The identity $f_\alpha(\tau^*(tm))=t f_\alpha(\tau^*m)$ shows that $e:=t\circ h-h\circ t_M$ is $R[\tau]$-linear. We will show that $e$ takes values in $R$ and corresponds to $[\alpha]_n$. We compute
\[
\text{for~}i\in \{1,\ldots,n-1\}:\quad e(\kappa_i)=th((t-\theta)^{i-1}\kappa_1)-h(t(t-\theta)^{i-1} \kappa_1)=0,
\]
and for $i=n$,
\begin{align*}
e(\kappa_n) &= th((t-\theta)^{n-1}\kappa_1)-h(t(t-\theta)^{n-1}\kappa_1) \\
&=-h(t(t-\theta)^{n-1}\kappa_1)=-h((t-\theta)^n\kappa_1)+h(\theta (t-\theta)^{n-1}\kappa_1) \\
&=-h(\tau \kappa_1)\\
&=-\alpha.
\end{align*}
Because $e$ is $R[\tau]$-linear and its values on the $R[\tau]$-basis $\kappa_1,\ldots,\kappa_n$ lie in $R$, it takes values in $R$. Hence $e(m)=m([\alpha]_n)$ for all $m\in M(C^{\otimes n})$, as desired.

Because $\Ext^{1,\eff}_{\AMot_R}(\underline{M}(C^{\otimes n}),\mathbbm{1})=\Ext^{1,\reg}_{\AMot_R}(\underline{M}(C^{\otimes n}),\mathbbm{1})$ is generated as an $\bF[t]$-module by the extensions $\iota(\alpha)$ for $\alpha\in R$ (see \cite[Proposition~5.3]{gazda}), we deduce that $C^{\otimes n}(R)$ is generated by the symbols $[\alpha]_n$. 
\end{Remark}

\section{Anderson generating functions}\label{sec:anderson-generating-functions}
Anderson used generating functions in \cite[\S 3]{anderson} to relate scattering matrices to uniformizability of $t$-modules. They seem to have acquired their name much later in the work of El-Guindy--Papanikolas \cite{elguindy-papanikolas}. The terminology ``generating'' evokes the construction of such series, which are built as generating series of $t$-division sequences of an Anderson $t$-module. In the setting of the Carlitz module, the archetype is the following series in~$t$
\[
f_C(z;t):=-\sum_{i=0}^{\infty}{\exp_C\left(\frac{z}{\theta^{i+1}}\right) t^i},
\]
where $z$ is a parameter in $\Lie_C(\CI)=\CI$. Using the $t$-linear exponential $\widehat{\exp}_C$, we can equivalently write
\begin{equation}\label{eq:generating-functions-exp}
f_C(z;t)=\widehat{\exp}_C\left(\frac{z}{t-\theta}\right).
\end{equation}
Let $\tilde{\pi}$ be a chosen Carlitz period. For $z=-\tilde{\pi}$, $f_C(z;t)$ recovers the Anderson--Thakur special function $\omega(t)$. 

In \cite{gazda-maurischat,gazda-maurischat2}, Maurischat and the author laid the foundations for a theory of special functions of an arbitrary Anderson module. In this section, we proceed likewise for generating functions, extending the results of \emph{loc.~cit.} We take the functional equation satisfied by generating functions as the defining property, resulting in an $A$-module of \emph{generating functions}
\[
\mathfrak{gf}_E(\CI)
\]
which contains special functions $\mathfrak{sf}_E(\CI)$ as a submodule. Given a parameter $z\in \Lie_E(\CI)$, we rather associate an element $f_E(z)$ of $\mathfrak{gf}_E(\CI)\otimes_A \Omega_{A/\bF}^1$ which, in the case of Carlitz, recovers $f_C(z;t)dt$.

The general theory of generating functions will enable us to offer a new perspective on the following theorem, originally due to Sinha \cite[Theorem~A.0.7]{sinha} in the $\bF[t]$-setting:
\begin{Theorem}\label{thm:sinha}
Let $B$ be a Banach $\FI$-algebra and let $E$ be a (not necessarily abelian) Anderson $A$-module over $\FI$, also denoted by $E$ after base change to $B$. There is a commutative diagram
\begin{equation}
\begin{tikzcd}
\Hom_{A\otimes B[\tau]}\left(M(E),\displaystyle\frac{B\langle A \rangle}{A\otimes B}\right)\otimes_A \Omega_{A/\bF}^1 \arrow[r,"\delta"]\arrow[d,"\wr"] & \Ext^{1}_{A\otimes B[\tau]}(M(E),A\otimes B)\otimes_A \Omega_{A/\bF}^1 \arrow[d,"\wr"',"\text{Thm~}\ref{thm:barsotti-weil}"] \\
\Lie_E(B) \arrow[r,"\exp_E"] & E(B)
\end{tikzcd}
\end{equation}
where the top horizontal arrow is the connecting homomorphism obtained after applying the functor $\Hom_{A\otimes B[\tau]}(M(E),-)$ to the short exact sequence
$0\to A\otimes B \to B\langle A \rangle \to \frac{B\langle A \rangle}{A\otimes B}\to 0$.
\end{Theorem}

Theorem \ref{thm:sinha} is obtained by showing that the residue at $\mathfrak{j}$ of $f_E(z)$ is $z$ while its ``residue at $\infty$'' is $-\exp_E(z)$ (we explain later what it means to take the residue at $\infty$ of such Tate series, hence the quotation marks). For the Carlitz module, this is recovered from the following formula borrowed from \cite{elguindy-papanikolas}, which holds for small $z$:
\begin{equation}\label{eq:elguindy-papanikolas-formula}
f_{C}(z;t)=\sum_{i=0}^{\infty}{\frac{\exp_C(z)^{q^i}}{(t-\theta)(t-\theta^q)\cdots (t-\theta^{q^i})}}.
\end{equation}
In our general setting, no comparably elegant formula seems to be available.

\subsection{Preliminary notation}
Throughout this section, let $B$ be a Banach $\FI$-algebra and let $E$ be an Anderson module over $\FI$; we use the same notation for its base change to $B$ and write $M(E)$ for the associated $B[\tau]$-module. Since $\FI$ is a finite product of fields, condition~\ref{item:zariski-form} gives an isomorphism $E\simeq\mathbb{G}_{a,\FI}^d$ of $\bF$-module schemes, which extends to an isomorphism over $B$. In particular, global coordinates over $B$ exist. Recall that the $A$-module $E(B)$ of $B$-points of $E$ carries a metric by identifying it with the group $\Hom_{B[\tau]}(M(E),B)$ and choosing a finite set of generators of $M(E)$ as a $B[\tau]$-module. This metric depends on the choice of the generators, hence is not canonical. However, the set of null sequences in $E(B)$ does not depend on such choices. Consequently, we may define
\[
E(B\langle A \rangle) := A\hat{\otimes} E(B)
\]
as the completion of $A\otimes E(B)$ with respect to the zero sequences in $E(B)$ (see \cite[\S 2]{gazda-maurischat}). This is consistent with the extended functor of points notation \eqref{eq:functor-of-points}, where $B\langle A\rangle$ is treated as a left $B[\tau]$-module where $\tau$ acts $A$-linearly:
\begin{equation}\label{eq:tate-valued-E}
E(B\langle A \rangle) \cong \Hom_{B[\tau]}(M(E),B\langle A \rangle)
\end{equation}
(see \cite[Theorem~3]{gazda-maurischat} and its proof). In the sequel, we will also consider $E(\sA)$ where $\sA$ is one of the following algebra of functions, treated as left $B[\tau]$-modules:
\begin{enumerate}[label=$\bullet$]
\item \emph{Classical Tate algebra.} Let $u\in A$ be nonconstant. We write $B\langle u\rangle$ for the Tate $B$-algebra in the variable $u$. It is a $B[\tau]$-submodule of $B\langle A \rangle$.
\item \emph{Rings of entire functions.} We write $\mathfrak{U}$ for the rigid analytification of $C\setminus \{\infty\}$ base changed to $B$ and we denote by $\mathcal{O}(\mathfrak{U})$ the algebra of analytic functions on $\mathfrak{U}$ (see \cite[\S 3.2]{gazda2} for a definition). This is canonically a left $B[\tau]$-module, where $\tau$ acts by $f\mapsto f^{(1)}$. Given $u$ non constant, it contains the subalgebra $B\langle\!\langle u\rangle\!\rangle$ of series in $u$ and coefficients in $B$ with infinite convergence radius. The following infinite product converges therein:
\[
\Pi_u:=\prod_{i\geq 0}{\left (1-\frac{u\otimes 1}{1\otimes u^{q^i}}\right)}.
\]
\item \emph{Module of meromorphic functions with poles of bounded order along $\fj$ and iterates.} We let $\mathcal{M}_J(\mathfrak{U})$ denote the algebra of meromorphic functions on $\mathfrak{U}$ whose poles are supported at the points $\mathfrak{j}$ and its iterates by the $A$-linear Frobenius, see \emph{loc.~cit.} This is once again a left $B[\tau]$-module.
\end{enumerate}

If $\sA$ is an $A\otimes B[\tau]$-module, we will often use that $E(\sA\otimes_A \Omega_{A/\bF}^1)$ is canonically isomorphic to $E(\sA)\otimes_A \Omega_{A/\bF}^1$. For the latter, we point out that the tensor product is taken with respect to the $A$-module structure on $\sA$ (and not the one on $E$).

\begin{Construction}[Residue at $\fj$ of $E$-valued functions]
Consider the $B$-linear map
\begin{equation}\label{eq:residue-at-j}
\res_{\fj}:\mathcal{M}_J(\mathfrak{U})\otimes_A \Omega_{A/\bF}^1\longrightarrow B
\end{equation}
which extracts the residue at $\fj$. Note that if a meromorphic differential in the source can be written as $\tau(\eta)$, then $\tau(\eta)$ is regular at $\fj$ as otherwise $\eta$ would have a pole at $\fj^{(-1)}$ which is excluded by the condition on the poles of $\eta$. In particular, $\res_{\fj}(\tau(\eta))=0$. We deduce that \eqref{eq:residue-at-j} promotes to a map of $B[\tau]$-modules, where $\tau$ acts by zero on its target $B$. Applying the extended functor of points $E$ defined in \eqref{eq:functor-of-points}, we obtain a map
\begin{equation}\label{eq:E(residue-at-j)}
\res_{\fj}:E(\mathcal{M}_J(\mathfrak{U}))\otimes_A \Omega_{A/\bF}^1 \longrightarrow \Lie_E(B).
\end{equation}
The map \eqref{eq:E(residue-at-j)} is the ``residue at $\fj$ of $E$-valued functions'' considered in \cite[\S 5.1.2]{gazda-maurischat2}.
\end{Construction}

\subsection{Analytic generating functions}
By analogy with the module of special functions defined in \cite{gazda-maurischat}, we introduce the following definition. 
\begin{Definition}[Analytic generating functions]
The module of \emph{Anderson generating functions of $E$ over $B$} is the $A$-module
\[
\mathfrak{gf}_E(B):=\{\Psi\in E(B\langle A \rangle)~|~\text{for~all~}a\in A:~(a\otimes 1-1\otimes \phi_a)\Psi\in A\otimes E(B)\}.
\]
\end{Definition}
\begin{Remark}
We distinguish two canonical submodules:
\begin{enumerate}[label=$\bullet$]
\item The module of \emph{special functions} $\mathfrak{sf}_E(B)$ of \cite{gazda-maurischat,gazda-maurischat2} is the submodule of elements $\omega\in E(B\langle A \rangle)$ such that $(a\otimes 1-1\otimes \phi_a)\omega=0$ for all $a\in A$.
\item $A\otimes E(B)$ is a submodule of $\mathfrak{gf}_E(B)$.
\end{enumerate}
\end{Remark}

Let $u\in A$ be a separating element and let $z\in \Lie_E(B)$. We consider the following expression
\[
f_E(z;u):=\widehat{\exp}_E\left((u\otimes 1-1\otimes \partial \phi_u)^{-1}z \right)\in E(B\langle u \rangle)
\]
which is similar to those originally considered by Anderson and generalizes the Carlitz-module expression \eqref{eq:generating-functions-exp}. The element $f_E(z;u)$ is well-defined because $u\otimes 1-1\otimes \ell(u)$ is invertible in $B\langle A\rangle$, \emph{e.g.} \cite[Lemma~3.15]{gazda}; hence $u\otimes 1-1\otimes \partial \phi_u$ is invertible as an endomorphism of $\Lie_E(B\langle A\rangle)$. In addition, it satisfies the relation 
\[
(u\otimes 1 -1\otimes \phi_u)f_E(z;u)=\exp_E(z),
\]
from which we deduce that $(a\otimes 1-1\otimes \phi_a)f_E(z;u)\in A\otimes E(B)$ for all $a\in \bF[u]$. The following series expansion is derived as in \cite[Lemma~5.7]{gazda-maurischat2}.
\begin{Lemma}\label{lem:converging-sum-for-psi}
Choose $N\geq 1$ such that $\mathfrak{j}^N\Lie_E(B)=0$, where $A$ acts on $\Lie_E(B)$ through $\partial\phi$. Such an $N$ exists because $A$ is finitely generated over $\bF$ and the corresponding nilpotent operators commute. Fix a choice of coordinates $\kappa:E\xrightarrow{\sim}\mathbb{G}^d_{a,B}$ and denote by $e_i^{\kappa}\in \mathrm{Mat}_{d}(B)$, $i\geq 0$, the coefficients of $\exp_E$ under $\kappa$. Write $\Lie_{\kappa}:\Lie_E(B) \to B^d$ for the induced map of Lie spaces. In the coordinates $\kappa$, we have
\begin{equation}\label{eq:expression-for-psiu}
f_E(z;u)=\sum_{j=0}^{N-1}\sum_{i\geq 0}{\frac{(1\otimes e_i^{\kappa})[(1\otimes \Lie_{\kappa})(1\otimes n_u^j)(z)]^{(i)}}{(u\otimes 1-1\otimes \ell(u)^{q^i})^{j+1}}}
\end{equation}
where $n_u:=\partial\phi_u-\ell(u)$ satisfies $n_u^N=0$. In particular, $f_E(z;u)$ belongs to $E(\Pi_u^{-N}B\langle\!\langle u \rangle\!\rangle)$.
\end{Lemma}

In general, the inclusion $\bF[u]\subset A$ is strict and $f_E(z;u)$ is not yet an element in $\mathfrak{gf}_E(B)$. However, appealing to the dualizing element of \eqref{eq:dualizing-element}, we can apply an averaging operation to derive one up to a twist by the module of K\"ahler differentials. This is the content of the next proposition. To state it, recall that we denoted by
\[
\sum_{i}{d_i\otimes_A \omega_i} \in \sD_{A/\bF[u]}\otimes_A \Omega_{A/\bF}^1
\]
the \emph{dualizing element} of \eqref{eq:dualizing-element}. For each $i$, choose a lift $\sum_j a_{ij}\otimes b_{ij}\in A\otimes A$ of $d_i$.

\begin{Proposition}[Constructing all generating functions]\label{prop:describe-generating-functions}
Consider the expression
\[
f_E(z):=\sum_{i,j}{a_{i,j}\omega_i \otimes \phi_{b_{i,j}} f_E(z;u)} \in E(\Omega_{A/\bF}^1\otimes_A B\langle A \rangle).
\]
Then,
\begin{enumerate}[label=$(\roman*)$]
\item\label{item:psi-generating} $f_E(z)$ belongs to $\mathfrak{gf}_E(B)\otimes_A \Omega_{A/\bF}^1$,
\item\label{item:iso} The map
\begin{equation}\label{eq:iso-generating-functions}
\Lie_E(B)\longrightarrow \frac{\mathfrak{gf}_E(B)}{A\otimes E(B)}\otimes_A \Omega_{A/\bF}^1, \quad z \longmapsto f_E(z)+\Omega_{A/\bF}^1\otimes E(B)
\end{equation}
is an $A$-linear isomorphism.
\item\label{item:reciprocal} The inverse of \eqref{eq:iso-generating-functions} is obtained by taking the residue at $\fj$. In particular, \eqref{eq:iso-generating-functions} does not depend on the choice of $u$.
\end{enumerate} 
\end{Proposition}

\begin{proof}[Proof of Proposition \ref{prop:describe-generating-functions}]
Consider the $\bF[u]$-linear map
\[
f_E(-;u):\Lie_E(B)\longrightarrow \frac{\{\psi\in E(B\langle u \rangle)~|~(u\otimes 1-1\otimes \phi_u) \psi\in \bF[u]\otimes_{\bF} E(B)\}}{\bF[u]\otimes E(B)}
\]
which assigns to $z\in \Lie_E(B)$ the expression $f_E(z;u)$ modulo $\bF[u]\otimes E(B)$. To show that this map is surjective, let $y$ be an element in the target. Because $\bF[u]\otimes E(B)$ is dense in $E(B\langle u \rangle)$, up to subtracting an element of $\bF[u]\otimes E(B)$ from $y$ we may assume that $y$ belongs to a ball of small enough radius in $E(B\langle u \rangle)$ on which the $\bF[u]$-linear exponential $\widehat{\exp}_E$ is an isometry. Assuming so, we can write $y=\widehat{\exp}_E(g)$ for some $g\in \Lie_E (B\langle u \rangle)$. Put $h:=(u\otimes 1-1\otimes\partial\phi_u)g=\sum_{n\geq 0}u^n\otimes h_n$, where $h_n\in\Lie_E(B)$ and $h_n\to 0$. Since $\widehat{\exp}_E(h)=(u\otimes 1-1\otimes\phi_u)y$ belongs to $\bF[u]\otimes E(B)$, we have $\exp_E(h_n)=0$ for all sufficiently large $n$. The exponential is injective on a sufficiently small neighborhood of zero, so $h_n=0$ for all sufficiently large $n$. Thus $h\in\bF[u]\otimes\Lie_E(B)$. We deduce that $g$ satisfies
\[
(u\otimes 1-1\otimes \partial \phi_u)g=\sum_{i}{p_i(u)\otimes z_i}, \quad \text{and} \quad y=\widehat{\exp}_E\left((u\otimes 1-1\otimes \partial \phi_u)^{-1} \sum_{i}{p_i(u)\otimes z_i}\right).
\]
Consequently,
\begin{align*}
y &=\widehat{\exp}_E\left((u\otimes 1-1\otimes \partial \phi_u)^{-1} \sum_{i}{p_i(u)\otimes z_i}\right) \\
&\equiv \widehat{\exp}_E\left((u\otimes 1-1\otimes \partial \phi_u)^{-1} \sum_{i}{1\otimes \partial \phi_{p_i(u)}(z_i)}\right) \pmod{\bF[u]\otimes E(B)} \\
&= f_E\left(\sum_i{\partial \phi_{p_i(u)}(z_i)};u\right) \pmod{\bF[u]\otimes E(B)}.
\end{align*}
This shows that $f_E(-;u)$ is surjective. It is also injective because one deduces from \eqref{eq:expression-for-psiu} that the residue at the point $u\otimes 1=1\otimes u$ of the differential $f_E(z;u)du$ is $z$. 

Extending scalars in $f_E(-;u)$ along $\bF[u]\to A$, we obtain an isomorphism
\begin{equation}\label{eq:Alinear}
A\otimes_{\bF[u]}\Lie_E(B)\xrightarrow{\sim} \frac{\{\psi\in E(B\langle A \rangle)~|~(u\otimes 1-1\otimes \phi_u) \psi\in A\otimes_{\bF} E(B)\}}{A\otimes E(B)}
\end{equation}
Applying the functor of $A$-invariants of Definition \ref{def:functor-of-A-invariants}, we get
\[
\sD_{A/\bF[u]}\otimes_A \Lie_E(B) \xrightarrow{\sim} \frac{\mathfrak{gf}_E(B)}{A\otimes E(B)}
\]
Tensoring over $A$ with $\Omega_{A/\bF}^1$ and using the isomorphism \eqref{eq:multiplication-different}, we get 
\[
\Lie_E(B) \xrightarrow{\sim} \frac{\mathfrak{gf}_E(B)}{A\otimes E(B)} \otimes_{A}\Omega_{A/\bF}^1.
\]
One checks by hand that the above maps $z$ to $f_E(z)+\Omega_{A/\bF}^1\otimes E(B)$, proving at once \ref{item:psi-generating} and \ref{item:iso}. To prove \ref{item:reciprocal}, we write
\begin{align*}
f_E(z) &= \sum_{i,j}{(a_{ij}\omega_i \otimes \phi_{b_{ij}})\widehat{\exp}_E\left((u\otimes 1-1\otimes \partial \phi_u)^{-1}z\right)} \\
&=\sum_{i,j}{(a_{ij}b_{ij}\omega_i \otimes 1)\widehat{\exp}_E\left((u\otimes 1-1\otimes \partial \phi_u)^{-1}z\right)} \\ &\qquad -\sum_{i,j}{(a_{ij}\omega_i \otimes 1)\widehat{\exp}_E\left((b_{ij}\otimes 1-1\otimes \partial \phi_{b_{ij}})(u\otimes 1-1\otimes \partial \phi_u)^{-1}z\right)} \\
&= \widehat{\exp}_E\left((u\otimes 1-1\otimes \partial \phi_u)^{-1}z\right) d(u\otimes 1)+(\text{some~expression~regular~at}~\fj).
\end{align*}
In particular, only the first term computes the residue at $\fj$. Because $u\in A$ is separating, $u\otimes 1-1\otimes u$ is a uniformizing parameter at $\fj$ (\emph{cf.} \cite[Lemma~2.3]{gazda-maurischat2}). In the formula \eqref{eq:expression-for-psiu} of Lemma \ref{lem:converging-sum-for-psi}, the only pair of indices that contributes to the residue is $(i,j)=(0,0)$ and we compute
\[
\res_{\fj}(f_E(z))=\res_{\fj}~\widehat{\exp}_E\left((u\otimes 1-1\otimes \partial \phi_u)^{-1}z\right) d(u\otimes 1)=\res_{\fj}~z \frac{d(u\otimes 1-1\otimes u)}{u\otimes 1 -1\otimes u}=z.
\]
\end{proof}

From Proposition \ref{prop:describe-generating-functions}, one may compute the residue at $\fj$ of generating functions. We now turn to the residue at $\infty$.
\begin{Construction}[residue at $\infty$]\label{cons:residue-at-infty}
We construct an $A$-linear map 
\begin{equation}\label{eq:exponential-from-generating-functions}
-\res_{\infty}:\frac{\mathfrak{gf}_E(B)}{A\otimes E(B)}\otimes_A \Omega_{A/\bF}^1 \longrightarrow E(B)
\end{equation}
that we will interpret later in Subsection \ref{subsec:formal-generating-functions} as taking (minus) ``the residue at $\infty$'' (hence the notation). We begin with the assignment $\psi\mapsto (u\otimes 1-1\otimes \phi_u)\psi$ which gives an $A\otimes A$-linear map
\[
\{\psi\in E(B\langle A \rangle)~|~(u\otimes 1-1\otimes \phi_u)\psi \in A\otimes E(B)\}\longrightarrow A\otimes E(B)
\]
mapping $f_E(z;u)$ to $\exp_E(z)$. In addition, it induces an $A\otimes A$-linear map:
\[
\frac{\{\psi\in E(B\langle A \rangle)~|~(u\otimes 1-1\otimes \phi_u)\psi \in A\otimes E(B)\}}{A\otimes E(B)}\to \frac{A\otimes E(B)}{(u\otimes 1-1\otimes \phi_u)(A\otimes E(B))}\cong A\otimes_{\bF[u]}E(B).
\]
Applying the functor of $A$-invariants and twisting by $\Omega_{A/\bF}^1$, we obtain a morphism as in \eqref{eq:exponential-from-generating-functions}. 
\end{Construction}

The following proposition, which is the key to Sinha's Theorem \ref{thm:sinha}, is immediate from the construction.
\begin{Proposition}\label{prop:exponential-from-generating-functions}
The morphism of Construction \ref{cons:residue-at-infty} is the unique $A$-linear map
\begin{equation}
\frac{\mathfrak{gf}_E(B)}{A\otimes E(B)}\otimes_A \Omega_{A/\bF}^1 \longrightarrow E(B)
\end{equation}
mapping $f_E(z)$ to $\exp_E(z)$. In particular, we have a commutative diagram:
\[\begin{tikzcd}
	& {\displaystyle\frac{\mathfrak{gf}_E(B)}{A\otimes E(B)}\otimes_A\Omega_{A/\bF}^1} & \\
	{\Lie_E(B)} && {E(B)}
	\arrow["{\res_{\mathfrak{j}}}"',"\sim", from=1-2, to=2-1]
	\arrow["{-\res_{\infty}}", from=1-2, to=2-3]
	\arrow["{\exp_E}", from=2-1, to=2-3]
\end{tikzcd}\]
\end{Proposition}

\subsection{Formal generating functions}\label{subsec:formal-generating-functions}
In this section, we justify why Construction \ref{cons:residue-at-infty} can be interpreted as taking the residue at $\infty$. Let $K_{\infty}$ denote the complete local field of $C$ at $\infty$; we denote by $\cO_{\infty}$ its valuation ring and by $\fm_{\infty}$ its maximal ideal. As in the theory of isocrystals, we consider the rings:
\[
\sA_{\infty}(B):=(\cO_{\infty}\otimes B)^{\wedge}_{\fm_{\infty}}, \quad \sB_{\infty}(B):=K_{\infty}\otimes_{\cO_{\infty}}\sA_{\infty}(B). 
\]
The ring $\sB_{\infty}(B)$ allows us to define the module of \emph{formal generating functions}. 
\begin{Definition}[Formal generating functions]
The module of formal generating functions of $E$ over $B$ is the $A$-module
\[
\mathfrak{GF}_E(B):=\{\psi \in E(\sB_{\infty}(B))~|~\text{for~all~}a\in A:~(a\otimes 1-1\otimes \phi_a)\psi \in A\otimes E(B)\}.
\]
\end{Definition}

For elements in the module of formal generating functions it makes sense to consider the residue at $\infty$. We consider the map $\res_{\infty}:E(\sB_{\infty}(B))\otimes_A \Omega_{A/\bF}^1\to E(B)$ induced the map $\sB_{\infty}(B)\otimes_{A}\Omega_{A/\bF}^1\to B$, which extracts the residue at $\infty$, after applying $\Hom_{B[\tau]}(M,-)$ (observe indeed that taking the residue at $\infty$ is a $B[\tau]$-linear operation).

\begin{Proposition}
The residue at $\infty$ induces an $A$-linear isomorphism
\begin{equation}\label{eq:residue-iso}
\frac{\mathfrak{GF}_E(B)}{A\otimes E(B)}\otimes_A \Omega_{A/\bF}^1 \xrightarrow{\sim} E(B)
\end{equation}
\end{Proposition}
\begin{proof}
We first build a map $\varepsilon$ as in \eqref{eq:exponential-from-generating-functions} using the same method as in Construction \ref{cons:residue-at-infty}, and then show that it coincides with taking (minus) the residue at $\infty$. The former step is similar; given a separating element $u\in A$, we consider the canonical map
\begin{equation}\label{eq:formal-generating-functions-iso}
u\otimes 1-1\otimes \phi_u:\{\psi \in E(\sB_{\infty}(B))~|~(u\otimes 1-1\otimes \phi_u)\psi \in A\otimes E(B)\}\longrightarrow A\otimes E(B).
\end{equation}
It is an isomorphism as, given $e\in A\otimes E(B)$, the following series converges in $E(\sB_{\infty}(B))$
\[
\xi_u(e):=\sum_{n\geq 0}{(u^{-(n+1)}\otimes \phi_u^{n})(e)}
\]
to an element that satisfies $(u\otimes 1-1\otimes \phi_u)\xi_u(e)=e$. It is easily verified that $e\mapsto \xi_u(e)$ is an inverse for \eqref{eq:formal-generating-functions-iso}. 

The map \eqref{eq:formal-generating-functions-iso} induces a morphism between quotient modules where on the source one mods out by $A\otimes E(B)$ and, on the target, by $(u\otimes 1-1\otimes \phi_u)(A\otimes E(B))$. We then twist by $\Omega_{A/\bF}^1$ and apply the functor of $A$-invariants to get the map $\varepsilon$.

Explicitly, the inverse of $\varepsilon$ maps $e$ to $\Xi(e):=\sum_{i,j}{(a_{ij}\omega_i\otimes \phi_{b_{ij}})\xi_u(e)}$. We compute its residue at $\infty$:
\[
\res_{\infty}(\Xi(e)) = \sum_{\substack{n\geq 0 \\ i,j}}{\res_{\infty}(a_{ij}\omega_i u^{-(n+1)})\otimes \phi_{b_{ij}u^n}}(e).
\]
Note that $\res_{\infty}(a_{ij}\omega_i u^{-(n+1)})\in \bF$ and vanishes for large $n\geq 0$. In particular, the above sum is finite and the following expression
\[
S:=\sum_{\substack{n\geq 0 \\ i,j}}{\res_{\infty}(a_{ij}\omega_i u^{-(n+1)})u^n b_{ij}}
\]
belongs to $A$. On the other hand, $S$ computes the residue at $\infty$ of the differential $\eta$ of Lemma \ref{lem:eta}, and hence $S=-1$ by \ref{item:residue-at-infty}. Consequently,
\[
\res_{\infty}(\Xi(e))=\phi_S(e)=-e.
\]
\end{proof}

The following lemma was invoked in the proof.
\begin{Lemma}\label{lem:eta}
For all $i$, fix a lift $\sum_{j}{a_{ij}\otimes b_{ij}}\in A\otimes A$ of $d_i\in \sD_{A/\bF[u]}\subset A\otimes_{\bF[u]} A$ and consider the meromorphic differential
\[
\eta:=\sum_{i,j}{\frac{a_{ij}\omega_i\otimes b_{ij}}{u\otimes 1-1\otimes u}}
\]
on the $K^s$-curve $C_{K^s}:=C\times_{\bF} K^s$ where $K^s$ denotes a separable closure of $K$. Then
\begin{enumerate}[label=$(\roman*)$]
\item\label{item:class-of-eta} The class of $\eta$ modulo $\Omega_{A/\bF}^1\otimes K^s$ does not depend on the lifts of $d_i$;
\item\label{item:residue-at-j} The residue of $\eta$ at $\fj$ is $1$;
\item\label{item:residue-at-j-sigma} $\eta$ is regular at $\fj^{\sigma}$ where $\sigma$ is any $\bF(u)$-linear embedding $\sigma:K\to K^s$ different from the inclusion;
\item\label{item:residue-at-infty} The residue of $\eta$ at $\infty$ is $-1$. 
\end{enumerate}
\end{Lemma}
\begin{proof}
Two lifts of $d_i$ differ by a sum of elements of the form $ap(u)\otimes b-a\otimes p(u)b=(a\otimes b)(p(u)\otimes 1-1\otimes p(u))$ for $a,b\in A$, $p(u)\in \bF[u]$. The difference is then a multiple of $u\otimes 1-1\otimes u$ which proves \ref{item:class-of-eta}. Note that, because $u$ is separating, the ideal $u\otimes 1-1\otimes u$ decomposes as $\prod_{\sigma}{\fj^{\sigma}}$ in the Dedekind domain $A\otimes K^s$, where the product is indexed over the $\bF(u)$-linear embeddings $\sigma:K\to K^s$ (see Lemma~2.4 and its proof in \cite{gazda-maurischat2}). In particular, $u\otimes 1-1\otimes u$ is a uniformizing parameter at $\fj$. By writing
\[
\eta=\sum_{i,j}{\frac{a_{ij}\omega_i b_{ij}\otimes 1}{u\otimes 1-1\otimes u}}-\sum_{i,j}{(a_{ij}\omega_i\otimes 1)\frac{b_{ij}\otimes 1-1\otimes b_{ij}}{u\otimes 1-1\otimes u}}
\]
we see that the former summand is $\frac{d(u\otimes 1)}{u\otimes 1-1\otimes u}$ and the latter is regular at $\fj$; because $d(u\otimes 1)=d(u\otimes 1-1\otimes u)$ in $\Omega_{A/\bF}^1\otimes K^s=\Omega_{A\otimes K^s/K^s}^1$, we get \ref{item:residue-at-j}. To prove \ref{item:residue-at-j-sigma}, we recall first that $\sD_{A/\bF[u]}$ is generated as an $A$-module by the elements $g_s\in A\otimes_{\bF[u]} A$ constructed as follows:
\begin{enumerate}
\item $s\in A$ is an element such that $K=\bF(u,s)$,
\item $f(X)\in \bF[u][X]$ is the minimal polynomial of $s$ over $\bF(u)$, 
\item $g(X,Y):=\frac{f(X)-f(Y)}{X-Y}\in \bF[u][X,Y]$,
\item $g_s:=g(s\otimes 1,1\otimes s)$.
\end{enumerate} 
Fix an $\bF(u)$-linear embedding $\sigma:K\to K^s$ which differs from the inclusion. As $f(s)=f(\sigma(s))=0$ but $s\neq \sigma(s)$, we see that $g_s$---and hence any element of $\sD_{A/\bF[u]}$---belongs to the kernel of the map $A\otimes_{\bF[u]} A\to K^s$, $a\otimes_{\bF[u]} b\mapsto a\sigma(b)$. By writing
\[
(u\otimes 1-1\otimes u)\eta=\sum_{i,j}{a_{ij}\omega_i\otimes b_{ij}}
\]
we deduce that the right-hand side vanishes modulo $\fj^{\sigma}$ whenever $\sigma$ differs from the inclusion. In particular $\eta$ is regular at $\fj^{\sigma}$. It remains to prove \ref{item:residue-at-infty}: it follows from \ref{item:residue-at-j-sigma} that $\fj$ and $\infty$ are the only possible poles of $\eta$. Using that the sum of the residues of a meromorphic differential on $C_{K^s}$ vanishes, we deduce it from \ref{item:residue-at-j}.
\end{proof}

\begin{Remark}
We finally explain why Construction \ref{cons:residue-at-infty} can be thought of as taking the residue at $\infty$. There is an abstract map
\begin{equation}\label{eq:formalization}
\frac{\mathfrak{gf}_E(B)}{A\otimes E(B)}\otimes_A \Omega_{A/\bF}^1 \longrightarrow \frac{\mathfrak{GF}_E(B)}{A\otimes E(B)}\otimes_A \Omega_{A/\bF}^1
\end{equation}
obtained by taking the negative of the composition of \eqref{eq:exponential-from-generating-functions} with the inverse of \eqref{eq:residue-iso}. Interpreting \eqref{eq:formalization} as a ``formalization map'' implies that \eqref{eq:exponential-from-generating-functions} should be thought of as ``minus \eqref{eq:residue-iso}'', hence minus taking the residue at $\infty$. 
\end{Remark}

\subsection{Exponential map as a connecting homomorphism}
We are now in a position to prove Sinha's theorem. We begin with a proposition.
\begin{Proposition}\label{prop:generating-functions-delta}
There is a commutative diagram
\[
\begin{tikzcd}
\displaystyle\frac{\mathfrak{gf}_E(B)}{A\otimes E(B)}\otimes_A \Omega_{A/\bF}^1 \arrow[r,"\sim"]\arrow[d,"-\res_{\infty}"'] &  \Hom_{A\otimes B[\tau]}\left(M(E),\displaystyle\frac{B\langle A \rangle}{A\otimes B}\right)\otimes_A \Omega_{A/\bF}^1\arrow[d,"\delta"] \\
E(B) \arrow[r,"\sim","\text{Thm}~\ref{thm:barsotti-weil}"'] & \Ext^{1}_{A\otimes B[\tau]}(M(E),A\otimes B)\otimes_A \Omega_{A/\bF}^1
\end{tikzcd}
\]
where $\delta$ is the connecting homomorphism obtained after applying the functor $\Hom_{A\otimes B[\tau]}(M(E),-)$ to the short exact sequence
\begin{equation}\label{eq:short-exact-sequence-tate}
0\longrightarrow A\otimes B \longrightarrow B\langle A \rangle \longrightarrow \frac{B\langle A \rangle}{A\otimes B}\longrightarrow 0.
\end{equation}
\end{Proposition}
\begin{proof}
We already mentioned that the identification $E(B)=\Hom_{B[\tau]}(M(E),B)$ promotes to $A\otimes E(B)=\Hom_{B[\tau]}(M(E),A\otimes B)$ and $E(B\langle A \rangle)=\Hom_{B[\tau]}(M(E),B\langle A \rangle)$ as $A\otimes A$-modules. Because $M(E)$ is projective over $B[\tau]$, we also get 
\begin{equation}\label{eq:quotient-of-identifications}
\frac{E(B\langle A \rangle)}{A\otimes E(B)} \xrightarrow{\sim} \Hom_{B[\tau]}\left(M(E),\frac{B\langle A \rangle}{A\otimes B}\right).
\end{equation}
Let $u\in A$ be a separating element. We identify on both sides of \eqref{eq:quotient-of-identifications} the submodule of elements that are killed by $u\otimes 1-1\otimes u$. This results in
\begin{equation}\label{eq:identification-hom}
\frac{\{\psi\in E(B\langle A \rangle)~|~(u\otimes 1-1\otimes \phi_u)\psi \in A\otimes E(B)\}}{A\otimes E(B)}\xrightarrow{\sim} \Hom_{B[u,\tau]}\left(M(E),\frac{B\langle A \rangle}{A\otimes B}\right).
\end{equation}
Multiplication by $u\otimes 1-1\otimes u$ defines a map from the source of \eqref{eq:identification-hom} to $ A\otimes_{\bF[u]} E(B)$ (identified with the cokernel of $u\otimes 1-1\otimes \phi_u$ on $A\otimes E(B)$) and from the target to 
\begin{equation}\label{eq:extensions-B[u,tau]}
\frac{\Hom_{B[\tau]}(M(E),A\otimes B)}{\{(u\otimes 1)g-g((1\otimes u)\cdot)~|~g\in \Hom_{B[\tau]}(M(E),A\otimes B)\}}.
\end{equation}
The module \eqref{eq:extensions-B[u,tau]} computes extensions between $M(E)$ and $A\otimes B$ in the category of $B[u,\tau]$\nobreakdash-modules. We obtain a commutative square
\begin{equation}\label{eq:diag-relation-to-barsotti-weil-1}
\begin{tikzcd}
\displaystyle\frac{\{\psi\in E(B\langle A \rangle)~|~(u\otimes 1-1\otimes \phi_u)\psi \in A\otimes E(B)\}}{A\otimes E(B)} \arrow[r,"\sim"] \arrow[d] & \Hom_{B[u,\tau]}\left(M(E),\frac{B\langle A \rangle}{A\otimes B}\right) \arrow[d] \\
A\otimes_{\bF[u]} E(B) \arrow[r] & \Ext^1_{B[u,\tau]}(M(E),A\otimes B)
\end{tikzcd}
\end{equation}
where the lower map is induced by the identification $A\otimes E(B)=\Hom_{B[\tau]}(M(E),A\otimes B)$, after passing to the corresponding quotients. To make $\mathfrak{gf}_E(B)$ appear in the upper left corner and rather get extensions of $A\otimes B[\tau]$-modules in the lower right corner, we apply the functor of $A$-invariants. We obtain the square of $A$-modules
\begin{equation}\label{eq:diag-relation-to-barsotti-weil-2}
\begin{tikzcd}
\displaystyle\frac{\mathfrak{gf}_E(B)}{A\otimes E(B)} \arrow[r,"\sim"] \arrow[d] & \Hom_{A\otimes B[\tau]}\left(M(E),\frac{B\langle A \rangle}{A\otimes B}\right) \arrow[d] \\
\sD_{A/\bF[u]}\otimes_A E(B) \arrow[r] & \Ext^1_{A\otimes B[\tau]}(M(E),A\otimes B)
\end{tikzcd}
\end{equation}
which is the claimed diagram after twisting \eqref{eq:diag-relation-to-barsotti-weil-2} by $\Omega_{A/\bF}^1$. It is routine from the proof of Theorem \ref{thm:barsotti-weil} that the lower horizontal map is the Barsotti--Weil isomorphism and that the right-hand vertical map is $\delta$. 
\end{proof}

\begin{proof}[Proof of Theorem \ref{thm:sinha}]
This follows immediately from Propositions \ref{prop:exponential-from-generating-functions} and \ref{prop:generating-functions-delta}.
\end{proof}

\section{Proof of the main theorem}\label{sec:proof}
Let $E$ be an abelian and uniformizable Anderson module over $\cO_F$ and $\underline{M}(E)$ be its motive. We write $\underline{M}=\underline{M}(E)^{\vee}$ for the dual of $\underline{M}(E)$ taken in the category of $A$-motives over $\cO_F$. This section is devoted to the proof of Theorem \ref{thm:comparison-taelman}.

\begin{proof}[Proof of Theorem \ref{thm:comparison-taelman}]
As before, we write $\underline{M}=\underline{M}(E)^{\vee}$ where $E$ is an abelian Anderson $A$\nobreakdash-module over $\cO_F$ which is uniformizable (over $F\otimes_K \CI$). Given a flat $\cO_F$-algebra $R$ and a $\KI$-Banach algebra $B$ which contains $R$, we consider the following diagram
\begin{equation}\label{eq:snake-diag}
\begin{tikzcd}
0 \arrow[r] & M\otimes_{\cO_F} R \arrow[r]\arrow[d,"\id-\tau_M"] & M\otimes_{A\otimes \cO_F} B\langle A \rangle \arrow[r]\arrow[d,"\id-\tau_M"] & M\otimes_{A\otimes \cO_F} \frac{B\langle A \rangle}{A\otimes R} \arrow[r]\arrow[d,"\id-\tau_M"] & 0  \\
0 \arrow[r] & N\otimes_{\cO_F} R \arrow[r] & N\otimes_{A\otimes \cO_F} B\langle A \rangle \arrow[r] & N\otimes_{A\otimes \cO_F} \frac{B\langle A \rangle}{A\otimes R} \arrow[r] & 0 
\end{tikzcd}
\end{equation}
where $N=M+\tau_M(\tau^*M)$. Note that both rows are exact; for the top one, this follows by tensoring the exact sequence of $A\otimes \cO_F$-modules
\begin{equation}\label{eq:short-exact-sequence-tate2}
0\longrightarrow A\otimes R \longrightarrow B\langle A \rangle \longrightarrow \frac{B\langle A \rangle}{A\otimes R} \longrightarrow 0
\end{equation}
with the flat module $M$. The bottom line requires an additional argument to show that $N\otimes_{\cO_F} R\to N\otimes_{A\otimes \cO_F} B\langle A \rangle$ is injective as $N$ may not be flat. Note that in $M[\fj^{-1}]$ we have an inclusion $N\subset \fj^{-n}M$ for $n$ large enough. Since $\cO_F\to R$ is flat, this further implies that $N\otimes_{\cO_F} R\to \fj^{-n}(M\otimes_{\cO_F} R)$ is injective; but the latter embeds in $M\otimes_{A\otimes \cO_F} B\langle A\rangle$ which agrees with $N\otimes_{A\otimes \cO_F} B\langle A \rangle$. This shows the exactness of the bottom row.

We denote by $\Snk(R;B)$ the long exact sequence obtained from the snake lemma applied to \eqref{eq:snake-diag}. Note that this is also the long exact sequence obtained by applying the functor $\Hom_{A\otimes \cO_F[\tau]}(M(E),-)$ to the short exact sequence \eqref{eq:short-exact-sequence-tate2}. In particular, applied with $R=B=\FI^s:=F\otimes_K \KI^s$, we get 
\[
\Snk(\FI^s;\FI^s)~:\quad 0 \longrightarrow \underline{M}_B \longrightarrow \Hom_{A\otimes \cO_F[\tau]}\left(M(E),\frac{\FI^s \langle A \rangle}{A\otimes \FI^s}\right) \xrightarrow{\delta} \Ext^{1,\reg}_{\AMot_{\FI^s}}(\mathbbm{1},\underline{M}_{\FI^s})\longrightarrow 0.
\]
Since $F_\infty^s$ is not itself a Banach algebra (but rather an ind-Banach algebra), the notation $F_\infty^s\langle A\rangle$ is understood as the direct limit
\[
F_\infty^s\langle A\rangle := \colim_{L/K_\infty~\text{finite~sep.}} (F\otimes_K L)\langle A\rangle
\]
(the sequence \eqref{eq:short-exact-sequence-tate2} with $R=B=\FI^s:=F\otimes_K \KI^s$ is exact, as a direct limit of exact sequences remains exact). By Theorem \ref{thm:sinha} --- applied to $B=F\otimes_K L$ and then taking the direct limit over finite separable extensions $L$ of $\KI$ --- we deduce that $\Snk(\FI^s;\FI^s)\otimes_A \Omega_{A/\bF}^1$ identifies with the exponential short exact sequence of $E$ over $\FI^s$.

On the other hand, the sequence $\Snk(\FI;\FI)$ has the following shape (\emph{cf.} \cite[Proposition~4.23]{gazda} and \cite[Theorem~3.16]{gazda2}):
\begin{equation}\label{eq:snakeFIFI}
0 \to \underline{M}_B^{G_{\infty}} \to \Hom_{A\otimes \cO_F[\tau]}\left(M(E),\frac{\FI \langle A \rangle}{A\otimes \FI}\right) \xrightarrow{\delta} \Ext^{1,\reg}_{\AMot_{\FI}}(\mathbbm{1},\underline{M}_{\FI})\to H^1(G_{\infty},\underline{M}_B)\to 0.
\end{equation}
There is no additional term on the right because $\id-\tau_M$ becomes surjective\footnote{Indeed, $N_{\FI}=N\otimes_{\cO_F} \FI$ is dense in $M\langle A \rangle := M\otimes_{A\otimes \cO_F} \FI\langle A \rangle=N\otimes_{A\otimes \cO_F} \FI\langle A \rangle$. In particular, given $y\in M\langle A \rangle$, we may write $y$ as $n+\varepsilon$ where $n\in N_{\FI}$ and $\varepsilon\in M\langle A \rangle$ sufficiently small that the series $x:=\sum_{k\geq 0}{\tau_M^k(\tau^{k*}\varepsilon)}$ converges in $M\langle A \rangle$. Then $(\id-\tau_M)x=\varepsilon\equiv y\pmod{N_{\FI}}$, so the class of $x$ maps to the class of $y$ in the quotient by $N_{\FI}$.} as a map from $M\otimes_{A\otimes \cO_F} \frac{B\langle A\rangle}{A\otimes B}$ to $$N\otimes_{A\otimes \cO_F} \frac{B\langle A\rangle}{A\otimes B}$$ for these choices of $R$ and $B$.

Since the functor of $G_{\infty}$-invariants is left-exact, the initial part of \eqref{eq:snakeFIFI} agrees with the left exact sequence obtained after applying $G_{\infty}$-invariants to $\Snk(\FI^s;\FI^s)$. In particular, $\Snk(\FI;\FI)\otimes_A \Omega_{A/\bF}^1$ identifies with the long exact sequence obtained by taking $G_{\infty}$-invariants of the exponential short exact sequence, and the map $\Ext^{1,\reg}_{\AMot_{\FI}}(\mathbbm{1},\underline{M}_{\FI})\to H^1(G_{\infty},\underline{M}_B)$ identifies, up to a twist by $\Omega_{A/\bF}^1$, with the connecting homomorphism \eqref{eq:connecting} obtained from the exponential short exact sequence. Finally, by \cite[Proposition~4.5]{gazda2}, $\Snk(\cO_F;\FI)$ has the form
\[
0 \to \underline{M}_B^{G_{\infty}} \to \Hom_{A\otimes \cO_F[\tau]}\left(M(E),\frac{\FI \langle A \rangle}{A\otimes \cO_F}\right) \to \Ext^{1,\reg}_{\AMot_{\cO_F}}(\mathbbm{1},\underline{M})\xrightarrow{r_B} H^1(G_{\infty},\underline{M}_B)\to \mathrm{Cl}(\underline{M})\to 0.
\]
We deduce from the functoriality of $\Snk$ that $r_B\otimes_A \mathrm{id}_{\Omega}$ identifies with $t_{\exp}$. This proves the theorem.
\end{proof}

\begin{proof}[Proof of Corollary \ref{cor:carlitz-symbols}]
Let $k_F$ be the full constant field of $F$ and let $g$ be its genus over $k_F$. For a place $v$ of $F$ above $\infty$, let $e_v$ be its ramification index, let $\deg_{k_F}(v)$ be its residue degree over $k_F$, and let $\operatorname{ord}_v$ be its integer-valued normalized valuation. Recall that $F_{\infty}=\prod_{v\mid\infty}F_v$ carries the norm $|x|:=\max_{v\mid\infty}q^{-\operatorname{ord}_v(x_v)/e_v}$. For $r>0$, let $B_r:=\{x\in F_\infty:|x|<r\}$ be the open ball of radius $r$. The Riemann--Roch theorem implies that if $c\in\bR$ satisfies
\begin{equation}\label{eq:inequality}
\sum_{v\mid\infty}(ce_v-1)\deg_{k_F}(v)>2g-2
\end{equation}
then $F_{\infty}=\cO_F+B_{q^c}$. In particular, if $n$ is such that \eqref{eq:inequality} holds with $c=\frac{qn}{q-1}$, then the same argument as in the proof of \cite[Proposition~2.2]{gazda-maurischat-ext} shows that the $\bF[t]$-module
\[
\Ext^{1,\reg}_{\cO_F,\infty}(\mathbbm{1},\underline{A}(n))
\]
is generated by the polylog classes $[\alpha]_n$ with $\alpha\in \cO_F$ such that $|\alpha|<q^{\frac{nq}{q-1}}$. We may conclude from Remark \ref{rem:compute-barsotti-weil-carlitz}.
\end{proof}

\section{The Hodge-Pink and Taelman regulators}\label{sec:regulators}
In this section, we show that the Taelman regulator of $E$ coincides with the Hodge-Pink regulator of $\underline{M}$ up to a twist by the module of K\"ahler differentials. Our technique is once again inspired by \cite{1-t-motif}, especially the proof of Theorem 8.1.2 in \emph{loc.~cit.}

Let $F$ be a finite extension of $K$ unramified at $\infty$ and denote by $\mathcal{O}_F$ the integral closure of $A$ in $F$. We write $F_{\infty}$ for the complete local field at the place $\infty$ seen in $F$. Let $E$ be an abelian and uniformizable Anderson module over $\cO_F$, $\Lambda_E$ its period lattice and $\Lambda_E^{+}$ its \emph{real} period lattice, \emph{i.e.} $\Lambda_E\cap \Lie_E(F_{\infty})$. Let us recall that \emph{the Taelman regulator of $E$} is the map
\[
\mathrm{Reg}_{T}(E): U(E/\mathcal{O}_F) \longrightarrow \Lie_E(F_{\infty})/(\KI\otimes_A\Lambda_E^{+})
\]
induced by the logarithm of $E$. Taelman's unit theorem implies that the image of $\mathrm{Reg}_{T}(E)$ is an $A$-lattice of full rank inside its target. For Drinfeld modules, his class formula relates its covolume to the special value of the $L$-series of $E$.

In \cite{gazda2}, the author introduced another regulator attached to an $A$-motive $\underline{M}$ over $F$, the \emph{Hodge-Pink regulator of $\underline{M}$}:
\[
\mathrm{Reg}_{HP}(\underline{M}):\Ext^{1,\reg}_{\cO_F,\infty}(\mathbbm{1},\underline{M}) \longrightarrow \Ext^{1,\mathrm{ha}}_{\infty}(\mathbbm{1},\mathcal{H}^+(\underline{M})).
\]
We show in this section that the latter generalizes the former.
 
\begin{Theorem}\label{thm:HP=T}
Let $\underline{M}=\underline{M}(E)^{\vee}_F$ be the base change to $F$ of the dual of the motive $\underline{M}(E)$ of $E$ taken in the category of $A$-motives over $\cO_F$. There is a commutative square of $A$-modules whose vertical maps are isomorphisms
\begin{equation}\label{eq:thm:HP=T}
\begin{tikzcd}[column sep=4em]
\Ext^{1,\reg}_{\cO_F,\infty}(\mathbbm{1},\underline{M})\otimes_A \Omega_{A/\bF}^1  \arrow[d,"\text{Thm. }\ref{thm:comparison-taelman}"']\arrow[r,"\mathrm{Reg}_{HP}(\underline{M})"] & \Ext^{1,\mathrm{ha}}_{\infty}(\mathbbm{1},\mathcal{H}^+(\underline{M}))\otimes_{A}\Omega_{A/\bF}^1 \arrow[d,"-\mathrm{res}_{\mathfrak{j}}"] \\
U(E/\mathcal{O}_F) \arrow[r,"\mathrm{Reg}_{T}(E)"] & \Lie_E(F_{\infty})/(\KI\otimes_A\Lambda_E^+)
\end{tikzcd}
\end{equation}
\end{Theorem}

Before beginning the proof of Theorem \ref{thm:HP=T}, we recall how the Hodge-Pink regulator is computed following the recipe given in \cite[Subsection~4.2]{gazda2}. Recall that \emph{the invariant Betti module of $\underline{M}$} is the $A$-module
\[
\underline{M}_B^+:=\{\omega\in M\otimes_{A\otimes F}\mathcal{M}_J(\mathfrak{U})~|~\tau_M(\tau^*\omega)=\omega\}
\]
where $\mathcal{M}_J(\mathfrak{U})$ denotes the ring of meromorphic functions on the analytification $\mathfrak{U}$ of $(C\setminus \{\infty\})_{F_{\infty}}$ whose poles are supported at $\mathfrak{j}$ and its $\tau$-iterates $\tau^{k*}\mathfrak{j}$, $k\geq 0$. The Laurent series expansion at $\mathfrak{j}$ of its elements lies in $(M+\tau_M(\tau^*M))\otimes_{A\otimes F}F_{\infty}[\![\mathfrak{j}]\!]$, interpreting $\underline{M}_B^+$ and its scalar extension $\underline{M}_B^+\otimes_A \KI$ as a submodule of the latter. It is proven in \cite[Theorem~3.32]{gazda2} that the extension space of Hodge-Pink structures is computed thanks to a natural isomorphism
\begin{equation}\label{eq:extension-hodge-pink}
\frac{(M+\tau_M(\tau^*M))\otimes_{A\otimes F}F_{\infty}[\![\mathfrak{j}]\!]}{\underline{M}_B^+\otimes_{A}\KI+M\otimes_{A\otimes F}F_{\infty}[\![\mathfrak{j}]\!]}\xrightarrow{\sim} \Ext^{1,\mathrm{ha}}_{\infty}(\mathbbm{1},\mathcal{H}^+(\underline{M})).
\end{equation}
For extensions of $A$-motives, write $M_{\cO_F}$ for $M(E)^{\vee}$ (as the notation suggests, this is the maximal integral model of $\underline{M}=\underline{M}(E)^{\vee}$ over $\cO_F$ from our good reduction hypothesis). As in \emph{loc.~cit.} set $N_{\mathcal{O}_F}:=M_{\mathcal{O}_F}[\fj^{-1}]\cap (M+\tau_M(\tau^*M))$ where the intersection is taken in $M[\mathfrak{j}^{-1}]$. We also consider the following natural extension of the invariant Betti module:
\[
\Delta_{\underline{M}}=\{\xi\in M\otimes_{A\otimes F} \mathcal{M}_J(\mathfrak{U})~|~\xi-\tau_M(\tau^* \xi)\in N_{\mathcal{O}_F}\}
\]
The Laurent series expansion at $\mathfrak{j}$ of the elements of $\Delta_{\underline{M}}$ lies in $(M+\tau_M(\tau^*M))\otimes_{A\otimes F}F_{\infty}[\![\mathfrak{j}]\!]$ as well, resulting in a map from the former to the latter:
\begin{equation}\label{eq:explicit-regulator}
\Delta_{\underline{M}}\longrightarrow (M+\tau_M(\tau^*M))\otimes_{A\otimes F}F_{\infty}[\![\mathfrak{j}]\!].
\end{equation}
Note that $\Delta_{\underline{M}}$ admits $\underline{M}_B^+$ and $M_{\cO_F}$ as submodules. In \cite[Subsection~4.2]{gazda2}, we derived a natural isomorphism
\begin{equation}\label{eq:extension-A-motives}
\frac{\Delta_{\underline{M}}}{\underline{M}_B^+ + M_{\mathcal{O}_F}} \xrightarrow{\sim} \Ext^{1,\reg}_{\cO_F,\infty}(\mathbbm{1},\underline{M}).
\end{equation}
It is proven in \cite[Proposition~4.8]{gazda2} that the Hodge-Pink regulator is given by minus \eqref{eq:explicit-regulator} under the descriptions \eqref{eq:extension-A-motives} and \eqref{eq:extension-hodge-pink}. 

The residue map in \eqref{eq:thm:HP=T} descends to the indicated quotient. Indeed, write $M_{\FI}[\![\fj]\!]:=M\otimes_{A\otimes F}\FI[\![\fj]\!]$ and $N_{\FI}[\![\fj]\!]:=(M+\tau_M(\tau^*M))\otimes_{A\otimes F}\FI[\![\fj]\!]$. The formal residue pairing identifies $(N_{\FI}[\![\fj]\!]/M_{\FI}[\![\fj]\!])\otimes_A\Omega_{A/\bF}^1$ with the $F_\infty$-linear dual of $M(E)_{F_\infty}/\tau_E(\tau^*M(E)_{F_\infty})$, hence with $\Lie_E(F_\infty)$. In particular, the residue vanishes on the regular term $M_{\FI}[\![\fj]\!]\otimes_A\Omega_{A/\bF}^1$. By the special-function case of Proposition~\ref{prop:describe-generating-functions}, the residue maps $\underline{M}_B^+\otimes_A\Omega_{A/\bF}^1$ to $\Lambda_E^+$. It therefore induces the stated map on quotients. The residue pairing is compatible with the $A$-actions and their continuous extensions to $\KI$, so this map is $\KI$-linear.

\begin{proof}[Proof of Theorem \ref{thm:HP=T}]
Let $x\in U(E/\mathcal{O}_F)$ be a Taelman unit. Then $x\in E(\mathcal{O}_F)$ and there exists $z\in \Lie_E(F_{\infty})$ such that 
\[
x=\exp_E(z).
\]
Taelman's regulator is induced by the logarithm of $E$, and is therefore given by 
\[
\mathrm{Reg}_T(x)= \log_E(x) \equiv z\pmod{\KI\otimes_A\Lambda_E^+}.
\]
Starting from $x$, we next follow the upper path in the square \eqref{eq:thm:HP=T}. We consider the Anderson generating function $f_E(z)\in \mathfrak{gf}_E(F_{\infty})\otimes_A\Omega_{A/\bF}^1$ of $z$ produced by Proposition \ref{prop:describe-generating-functions} inside the quotient module
\[
\frac{\mathfrak{gf}_E(F_{\infty})}{A\otimes E(F_{\infty})}\otimes_A \Omega_{A/\bF}^1 = \Hom_{A\otimes F_{\infty}[\tau]}\left(M(E),\frac{\mathcal{M}_J(\mathfrak{U})}{A\otimes F_{\infty}}\right)\otimes_A \Omega_{A/\bF}^1
\]
(the above identification with a homomorphism module follows from Proposition \ref{prop:generating-functions-delta} and the meromorphic continuation supplied by the construction of $f_E(z)$). Indeed, these classes exhaust the source by Proposition \ref{prop:describe-generating-functions}, and restriction to the Tate disk is injective. In particular, $f_E(z)\pmod{\Omega_{A/\bF}^1\otimes E(F_{\infty})}$ seen as a homomorphism can be lifted---because $E$ is abelian---to an $A\otimes F_{\infty}$-linear map
\[
(\xi_0:M(E)\longrightarrow \mathcal{M}_J(\mathfrak{U})\otimes_A\Omega_{A/\bF}^1) \in M\otimes_{A\otimes F} \mathcal{M}_J(\mathfrak{U}) \otimes_A \Omega_{A/\bF}^1. 
\]
By construction, we have that $\xi\circ \tau_E-\tau_{\mathbbm{1}}\circ \tau^*\xi$ factors as an $A\otimes F_{\infty}$-linear map from $\tau^*M(E)$ to $\Omega_{A/\bF}^1 \otimes F_{\infty}$. Because $\tau_M(f)=\tau_{\mathbbm{1}}\circ \tau^* f \circ \tau_E^{-1}$, it follows that $\xi$ belongs to $\Delta_{\underline{M}}\otimes_A \Omega_{A/\bF}^1$. Under \eqref{eq:extension-A-motives}, it represents the class of the extension associated to $x$ via the Barsotti--Weil isomorphism. In addition, the residue at $\fj$ of $\xi$ agrees with that of $f_E(z)$ and in particular is $z$ by Proposition~\ref{prop:describe-generating-functions}. Consequently, the composition
\[
\Delta_{\underline{M}}\otimes_A\Omega_{A/\bF}^1 \xrightarrow{-\eqref{eq:explicit-regulator}} \left(\frac{(M+\tau_M(\tau^*M))\otimes_{A\otimes F}F_{\infty}[\![\mathfrak{j}]\!]}{\underline{M}_B^+\otimes_A \KI+M\otimes_{A\otimes F}F_{\infty}[\![\mathfrak{j}]\!]}\right)\otimes_A\Omega_{A/\bF}^1 \xrightarrow{-\mathrm{res}_{\mathfrak{j}}} \Lie_E(F_{\infty})/\KI\otimes_A \Lambda_E^+
\]
maps $\xi$ to $z \pmod{\KI\otimes_A \Lambda_E^+}$, as desired.

This shows that the square \eqref{eq:thm:HP=T} commutes. It remains to prove that $-\mathrm{res}_{\fj}$ is an isomorphism. Extend the two horizontal maps $\KI$-linearly on their sources, keeping their targets, which are already $\KI$-vector spaces. Both horizontal maps then become isomorphisms: for the upper one, this is \cite[Theorem~6.3]{gazda2}, and for the lower one this follows, \emph{e.g.}, from \cite[Theorem~1.7]{fang}. The left vertical map remains an isomorphism after tensoring with $\KI$. Since the right vertical map is $\KI$-linear and the square commutes, it is an isomorphism as well.
\end{proof}


{\footnotesize
\begin{thebibliography}{GKK25}

\bibitem[And86]{anderson}
G.~W. Anderson,
\newblock $t$-motives,
\newblock \emph{Duke Math. J.} \textbf{53} (1986), no.~2, 457--502.
\newblock
\href{https://doi.org/10.1215/S0012-7094-86-05328-7}
{doi:10.1215/S0012-7094-86-05328-7}.

\bibitem[BVK16]{barbieri-viale-kahn2016}
L.~Barbieri-Viale and B.~Kahn,
\newblock On the derived category of $1$-motives,
\newblock \emph{Ast\'erisque}, no.~381, Soci\'et\'e Math\'ematique de
France, Paris, 2016, xi+254~pp.
\newblock
\href{https://doi.org/10.24033/ast.998}
{doi:10.24033/ast.998}.

\bibitem[BH07]{bockle-hartl}
G.~B{\"o}ckle and U.~Hartl,
\newblock Uniformizable families of $t$-motives,
\newblock \emph{Trans. Amer. Math. Soc.} \textbf{359} (2007), no.~8,
3933--3972.
\newblock
\href{https://doi.org/10.1090/S0002-9947-07-04136-0}
{doi:10.1090/S0002-9947-07-04136-0}.

\bibitem[CGL26]{caruso-gazda-lucas}
X.~Caruso, Q.~Gazda, A.~Lucas,
\newblock Class formulas over $A[z]$,
\newblock unpublished preprint,
\newblock \href{https://arxiv.org/abs/2609.17321}{arxiv:2609.17321}

\bibitem[Del74]{deligne1974}
P.~Deligne,
\newblock Th\'eorie de Hodge III,
\newblock \emph{Publ. Math. Inst. Hautes \'Etudes Sci.} \textbf{44}
(1974), 5--77.
\newblock
\href{https://doi.org/10.1007/BF02685881}
{doi:10.1007/BF02685881}.

\bibitem[EP14]{elguindy-papanikolas}
A.~El-Guindy and M.~A. Papanikolas,
\newblock Identities for Anderson generating functions for Drinfeld modules,
\newblock \emph{Monatsh. Math.} \textbf{173} (2014), no.~4, 471--493.
\newblock
\href{https://doi.org/10.1007/s00605-013-0543-9}
{doi:10.1007/s00605-013-0543-9}.

\bibitem[Fa15]{fang}
J. Fang,
\newblock Special \(L\)-values of abelian \(t\)-modules,
\newblock \emph{J. Number Theory} \textbf{147} (2015), 300--325.

\bibitem[Gaz22]{gazda2}
Q.~Gazda,
\newblock Regulators in the arithmetic of function fields,
\newblock unpublished preprint,
\href{https://arxiv.org/abs/2207.03461}
{arXiv:2207.03461} (2022; version~3, 2026).

\bibitem[Gaz24]{gazda}
Q.~Gazda,
\newblock On the integral part of $A$-motivic cohomology,
\newblock \emph{Compos. Math.} \textbf{160} (2024), no.~8, 1715--1783.
\newblock
\href{https://doi.org/10.1112/S0010437X24007218}
{doi:10.1112/S0010437X24007218}.

\bibitem[GM21]{gazda-maurischat}
Q.~Gazda and A.~Maurischat,
\newblock Special functions and Gauss--Thakur sums in higher rank and
dimension,
\newblock \emph{J. Reine Angew. Math.} \textbf{773} (2021), 231--261.
\newblock
\href{https://doi.org/10.1515/crelle-2020-0038}
{doi:10.1515/crelle-2020-0038}.

\bibitem[GM22]{gazda-maurischat-ext}
Q.~Gazda and A.~Maurischat,
\newblock Carlitz twists: their motivic cohomology, regulators, zeta values
and polylogarithms,
\newblock unpublished preprint,
\href{https://arxiv.org/abs/2212.02972}
{arXiv:2212.02972} (2022; version~2, 2023).

\bibitem[GM24]{gazda-maurischat2}
Q.~Gazda and A.~Maurischat,
\newblock Residue of special functions of Anderson $A$-modules at the
characteristic graph,
\newblock \emph{J. Number Theory} \textbf{260} (2024), 1--28.
\newblock
\href{https://doi.org/10.1016/j.jnt.2024.01.013}
{doi:10.1016/j.jnt.2024.01.013}.

\bibitem[GM25]{gazda-maurischat-pairing}
Q.~Gazda and A.~Maurischat,
\newblock Pairing Anderson motives via formal residues in the Frobenius
endomorphism,
\newblock unpublished preprint,
\href{https://arxiv.org/abs/2504.01926}
{arXiv:2504.01926} (2025; version~2, 2025).

\bibitem[GKK25]{gkk}
F.~G{\l}och, D.~E. Kedzierski, and P.~Kraso{\'n},
\newblock Algorithms for determination of $t$-module structures on some
extension groups,
\newblock \emph{Int. J. Number Theory} \textbf{21} (2025), no.~8,
1889--1922.
\newblock
\href{https://doi.org/10.1142/S1793042125500927}
{doi:10.1142/S1793042125500927}.
\newblock Corrigendum, published online 4 June 2026:
\href{https://doi.org/10.1142/S1793042126920028}
{doi:10.1142/S1793042126920028}.

\bibitem[Har19]{hartl-isogeny}
U.~Hartl,
\newblock Isogenies of abelian Anderson $A$-modules and $A$-motives,
\newblock \emph{Ann. Sc. Norm. Super. Pisa Cl. Sci. (5)}
\textbf{19} (2019), no.~4, 1429--1470.
\newblock
\href{https://doi.org/10.2422/2036-2145.201612_003}
{doi:10.2422/2036-2145.201612\_003}.

\bibitem[KK26]{kk}
D.~E. Kedzierski and P.~Kraso{\'n},
\newblock Weil--Barsotti formula for $t$-modules,
\newblock \emph{J. Number Theory} \textbf{278} (2026), 1--25.
\newblock
\href{https://doi.org/10.1016/j.jnt.2025.04.013}
{doi:10.1016/j.jnt.2025.04.013}.

\bibitem[KP24]{kp}
D.~E. Kedzierski and P.~Kraso{\'n},
\newblock On $\operatorname{Ext}^1$ for Drinfeld modules,
\newblock \emph{J. Number Theory} \textbf{256} (2024), 97--135.
\newblock
\href{https://doi.org/10.1016/j.jnt.2023.09.007}
{doi:10.1016/j.jnt.2023.09.007}.

\bibitem[Mor18]{mornev-shtuka}
M.~Mornev,
\newblock Shtuka cohomology and special values of Goss $L$-functions,
\newblock Ph.D. thesis, Universiteit Leiden, 2018.
\newblock
\href{https://hdl.handle.net/1887/61145}{hdl:1887/61145};
updated version:
\href{https://arxiv.org/abs/1808.00839}{arXiv:1808.00839}.

\bibitem[PR03]{PR}
M.~A. Papanikolas and N.~Ramachandran,
\newblock A Weil--Barsotti formula for Drinfeld modules,
\newblock \emph{J. Number Theory} \textbf{98} (2003), no.~2, 407--431.
\newblock
\href{https://doi.org/10.1016/S0022-314X(02)00047-1}
{doi:10.1016/S0022-314X(02)00047-1}.

\bibitem[Poo95]{poonen1995}
B.~Poonen,
\newblock Local height functions and the Mordell--Weil theorem for Drinfeld
modules,
\newblock \emph{Compos. Math.} \textbf{97} (1995), no.~3, 349--368.
\newblock
\href{https://www.numdam.org/item/CM_1995__97_3_349_0/}
{NUMDAM: CM\_1995\_\_97\_3\_349\_0}.

\bibitem[Sin97]{sinha}
S.~K. Sinha,
\newblock Periods of $t$-motives and transcendence,
\newblock \emph{Duke Math. J.} \textbf{88} (1997), no.~3, 465--535.
\newblock
\href{https://doi.org/10.1215/S0012-7094-97-08820-7}
{doi:10.1215/S0012-7094-97-08820-7}.

\bibitem[Tae20]{1-t-motif}
L.~Taelman,
\newblock $1$-$t$-motifs,
\newblock in \emph{$t$-Motives: Hodge Structures, Transcendence, and Other
Motivic Aspects} (G.~B{\"o}ckle, D.~Goss, U.~Hartl, and M.~Papanikolas,
eds.), EMS Ser. Congr. Rep., EMS Publishing House, Z{\"u}rich, 2020,
pp.~417--439.
\newblock
\href{https://doi.org/10.4171/198-1/8}
{doi:10.4171/198-1/8}.

\bibitem[Tae10]{taelman-dirichlet}
L.~Taelman,
\newblock A Dirichlet unit theorem for Drinfeld modules,
\newblock \emph{Math. Ann.} \textbf{348} (2010), no.~4, 899--907.
\newblock
\href{https://doi.org/10.1007/s00208-010-0506-6}
{doi:10.1007/s00208-010-0506-6}.

\bibitem[Woo95]{woo}
S.~S. Woo,
\newblock Extensions of Drinfeld modules of rank $2$ by the Carlitz module,
\newblock \emph{Bull. Korean Math. Soc.} \textbf{32} (1995), no.~2,
251--257.

\end{thebibliography}
}
\end{document}